\documentclass[thmsa,12pt,sw20lart]{article}%
\usepackage{amsmath}
\usepackage{amsfonts}
\usepackage{amssymb}
\usepackage{graphicx}%
\newtheorem{theorem}{Theorem}

\newtheorem{corollary}[theorem]{Corollary}

\newtheorem{lemma}[theorem]{Lemma}

\newtheorem{remark}[theorem]{Remark}

\newenvironment{proof}[1][Proof]{\noindent\textbf{#1.} }{\ \rule{0.5em}{0.5em}}
\begin{document}

\author{Yu-Chu Lin\ \&\ Dong-Ho Tsai\thanks{Mathematics\ Subject
Classification:\ 35K15, 35K55.}\thanks{Research supported by NSC\ (grant
number 95-2115-M-007-009)\ and\ National Center of Theoretical Sciences\ of
Taiwan.}}
\title{Nonlocal Flow of Convex Plane Curves\ and\ Isoperimetric Inequalities}
\maketitle

\begin{abstract}
In the first part of the paper we survey some nonlocal flows of convex plane
curves ever studied so far\ and discuss properties of the flows related to
enclosed\ area and length, especially the isoperimetric ratio and the
isoperimetric difference.\ We also study a new nonlocal flow\ of convex plane
curves and discuss its evolution behavior.\ 

In the second part of the paper we discuss necessary and sufficient conditions
(in terms of the (mixed)\ isoperimetric ratio or (mixed)\ isoperimetric
difference)\ for two convex closed curves to be homothetic\ or parallel.

\end{abstract}

\section{Introduction}

Recently there has been some interest in the nonlocal flow of \emph{convex}
closed plane curves. See the papers by Gage\ \cite{GA2},\ Jiang-Pan\ \cite{JP}%
,\ Pan-Yang \cite{PY},\ Ma-Cheng\ \cite{MC}\ and Ma-Zhu\ \cite{MZ}).\ All of
the above papers deal with the evolution of a given convex\footnote{In this
paper, "convex" always means "strictly convex".\ A convex closed
plane\ curve\ has\ positive curvature everywhere.\ Also, for simplicity, all
curves considered in this paper are smooth with positive orientation.}
embedded\ closed plane\ curve $\gamma_{0}$. The general form of the equation
is given by%
\begin{equation}
\left\{
\begin{array}
[c]{l}%
\dfrac{\partial X}{\partial t}\left(  \varphi,t\right)  =\left[  F\left(
k\left(  \varphi,t\right)  \right)  -\lambda\left(  t\right)  \right]
\mathbf{N}_{in}\left(  \varphi,t\right)
\vspace{3mm}%
\\
X\left(  \varphi,0\right)  =X_{0}\left(  \varphi\right)  ,\ \ \ \varphi\in
S^{1},
\end{array}
\right.  \label{main}%
\end{equation}
which is a parabolic initial value problem.\ Here $X_{0}\left(  \varphi
\right)  :S^{1}\rightarrow\gamma_{0}\ $is a smooth parametrization of
$\gamma_{0}$; $k\left(  \varphi,t\right)  \ $is the curvature of the evolving
curve $\gamma_{t}=\gamma\left(  \cdot,t\right)  \ $(parametrized by $X\left(
\varphi,t\right)  $)\ at the point $\varphi;\ $and $\mathbf{N}_{in}\left(
\varphi,t\right)  $ is the inward normal of $\gamma_{t}$. As for the
speed,$\ F\left(  k\right)  \ $is a given function\ of curvature\ satisfying
the\emph{ parabolic condition} $F^{\prime}\left(  z\right)  >0\ $for all
$z\ $in its\ domain$\ $and $\lambda\left(  t\right)  \ $is a function of time,
which may depend on certain global quantities of $\gamma_{t},\ $say length
$L\left(  t\right)  ,$ enclosed area $A\left(  t\right)  ,\ $or others\ (see
(\ref{F2-1})\ and\ (\ref{F4})).$\ $Note that if $\lambda\left(  t\right)  $
depends on $\gamma_{t},\ $then it is not known beforehand.\ 

The results claimed in each of the above mentioned papers are more or less the
same:\ the flows preserve the convexity of a given initial curve $\gamma_{0}$
and evolve it to a round circle (or round point)\ in $C^{\infty}$ sense as
$t\rightarrow\infty.\ $

For $k$-type flows,\ the following three flows:
\begin{equation}
F\left(  k\right)  -\lambda\left(  t\right)  =k-\frac{2\pi}{L\left(  t\right)
}\ \ \ \text{(area-preserving,\ gradient flow of\ }L^{2}-4\pi A\text{)}
\label{F1}%
\end{equation}
and
\begin{equation}
F\left(  k\right)  -\lambda\left(  t\right)  =k-\frac{L\left(  t\right)
}{2A\left(  t\right)  }\ \ \ \text{(gradient flow of }L^{2}/4\pi A\text{)}
\label{F2}%
\end{equation}
and \
\begin{equation}
F\left(  k\right)  -\lambda\left(  t\right)  =k-\frac{1}{2\pi}\int
_{0}^{L\left(  t\right)  }k^{2}ds\ \ \ \text{(length-preserving)} \label{F2-1}%
\end{equation}
are studied by\ Gage \cite{GA2}, Jiang-Pan\ \cite{JP},\ and\ Ma-Zhu\ \cite{MZ} respectively.\ 

For $1/k$-type\ flows, the following two flows:%
\begin{equation}
F\left(  k\right)  -\lambda\left(  t\right)  =\frac{1}{L\left(  t\right)
}\int_{0}^{L\left(  t\right)  }\frac{1}{k}ds-\frac{1}{k}%
\ \ \ \text{(area-preserving)} \label{F4}%
\end{equation}
and
\begin{equation}
F\left(  k\right)  -\lambda\left(  t\right)  =\frac{L\left(  t\right)  }{2\pi
}-\frac{1}{k}\ \ \ \text{(length-preserving)} \label{F3}%
\end{equation}
are studied by\ Ma-Cheng\ \cite{MC} and\ Pan-Yang \cite{PY} respectively.\ 

We know that (see \cite{GA2}) for a family of time-dependent
simple\ closed\ curves $X\left(  \varphi,t\right)  :S^{1}\times\lbrack
0,T)\rightarrow\mathbb{R}^{2}$ with \textsf{time variation}%
\begin{equation}
\frac{\partial X}{\partial t}\left(  \varphi,t\right)  =W\left(
\varphi,t\right)  \in\mathbb{R}^{2}, \label{W}%
\end{equation}
its length $L\left(  t\right)  $ and enclosed area $A\left(  t\right)
\ $satisfy the following:
\begin{equation}
\frac{dL}{dt}\left(  t\right)  =-\int_{\gamma_{t}}\left\langle W,\ k\mathbf{N}%
_{in}\right\rangle ds,\ \ \text{\ \ \ }\frac{dA}{dt}\left(  t\right)
=-\int_{\gamma_{t}}\left\langle W,\ \mathbf{N}_{in}\right\rangle ds.
\label{LA-evo}%
\end{equation}
This says that the well-known \emph{curve shortening flow} (with $F\left(
k\right)  =k,\ \lambda=0\ $in (\ref{main})) is the gradient flow of the length
functional.\ See Gage-Hamilton\ \cite{GH}). Also the \emph{unit-speed inward
normal flow} (with $F\left(  k\right)  =0,\ \lambda=-1\ $in (\ref{main})) is
the gradient flow of the area functional.

In particular, for convex plane curves evolution,\ we can check that flows
(\ref{F1}), (\ref{F4}) are \emph{area-preserving }and\ flows (\ref{F2-1}%
),\ (\ref{F3}) are \emph{length-preserving.}

As for flow (\ref{F2}), it$\ $is \emph{length-decreasing} due to Gage's
inequality for convex closed curves\ (see \cite{GA3}): \emph{ }
\begin{equation}
\frac{dL}{dt}\left(  t\right)  =-\int_{\gamma_{t}}k^{2}ds+\frac{\pi L\left(
t\right)  }{A\left(  t\right)  }\leq0. \label{G-ineq}%
\end{equation}
It is also\ \emph{area-increasing }due$\ $to (note that $%
{\displaystyle\int_{\gamma_{t}}}
kds=2\pi$):\emph{ }%
\[
\frac{dA}{dt}\left(  t\right)  =-2\pi+\frac{L^{2}\left(  t\right)  }{2A\left(
t\right)  }\geq0.
\]
Thus this flow is the most efficient in evolving a convex\ curve to a round
circle.\ As we shall see in Lemma \ref{lem0} below,\ Jiang-Pan's
flow\ (\ref{F2}) is \emph{the} \emph{gradient flow of the isoperimetric ratio
functional}.

As a comparison, for$\ k$-type flows\ with speed of the form$\ \left[
k-p\left(  t\right)  \right]  \mathbf{N}_{in},\ $we have
\begin{equation}
\dfrac{dL}{dt}\left(  t\right)  =-\int_{\gamma_{t}}k^{2}ds+2\pi p\left(
t\right)  ,\ \ \ \ \ \dfrac{dA}{dt}\left(  t\right)  =-2\pi+p\left(  t\right)
L \label{D1}%
\end{equation}
and for$\ 1/k$-type flows with speed of the form$\ \left[  q\left(  t\right)
-1/k\right]  \mathbf{N}_{in},\ $we have \
\begin{equation}
\dfrac{dL}{dt}\left(  t\right)  =-2\pi q\left(  t\right)  +L,\ \ \ \ \ \dfrac
{dA}{dt}\left(  t\right)  =-q\left(  t\right)  L+\int_{\gamma_{t}}\dfrac{1}%
{k}ds. \label{D2}%
\end{equation}
It is interesting to observe that when $q\left(  t\right)  =1/p\left(
t\right)  ,\ $there is a "\emph{dual relation}"\ between (\ref{D1})\ and
(\ref{D2}), i.e.,\
\begin{equation}
\frac{1}{q\left(  t\right)  }\dfrac{dL}{dt}\left(  t\right)  \ \text{(for
}1/k\text{-type\ flows)}=\dfrac{dA}{dt}\left(  t\right)  \ \text{(for
}k\text{-type\ flows).} \label{D3}%
\end{equation}
Hence in the above, flows\ (\ref{F1})\ and (\ref{F3})\ are dual.\ We shall
consider the dual flow of (\ref{F2}) in Section \ref{sec}.

In the first part of the paper, we observe some interesting behavior of a
general nonlocal flow\ (\ref{main}), especially the properties related to the
isoperimetric difference $L^{2}-4\pi A$ and isoperimetric ratio $L^{2}/4\pi
A.\ $We also discuss certain difficulty in dealing with the flow
(\ref{F3})\ and (\ref{F4}), especially the possibility of curvature blowing up
in finite time.\ 

In the second part, we discuss certain necessary and sufficient conditions (in
terms of the \emph{mixed isoperimetric\ ratio }$L_{1}L_{2}/4\pi A_{12}%
\ $and\emph{\ mixed isoperimetric difference\ }$L_{1}L_{2}-4\pi A_{12}$)\ for
two convex closed plane\ curves $\gamma_{1},\ \gamma_{2}\ $to be
\emph{homothetic }or\ \emph{parallel}.\ Here $A_{12}\ $is the \emph{mixed
area} determined by $A_{1},\ A_{2}.\ $

For simplicity, throughout the rest of the paper, we shall use the
following\ two abbreviations:%
\begin{equation}
\text{IPR}=\text{isoperimetric\ ratio,\ \ \ IPD}%
=\text{isoperimetric\ difference.}%
\end{equation}

\section{The decreasing of the IPD}

In this section we first\ prove an interesting property of the flow
(\ref{main}).\ It says that the IPD$\ L^{2}-4\pi A$ is always
non-increasing.\ To explain this, we need the following nice inequality due to
Andrews\ (see\ p. 341 of \cite{A}).\ One can view it as a generalization of
the classical H\"{o}lder inequality.\ 

\begin{lemma}
(\emph{Andrews's inequality})\ Let $M$ be a compact Riemannian manifold with a
volume form $d\mu,$\ and let $\xi$\ be a continuous function on$\;M.$\ Then
for any\textbf{ increasing} continuous function $F:\mathbb{R\rightarrow R}%
,\;$we have
\begin{equation}
\int_{M}\xi d\mu\int_{M}F\left(  \xi\right)  d\mu\leq\int_{M}d\mu\int_{M}\xi
F\left(  \xi\right)  d\mu. \label{Ben}%
\end{equation}
If $F$\ is strictly increasing, then equality holds if and only if $\xi$\ is a
constant function on $M.$\ Similarly,\ if $F:\mathbb{R}\rightarrow
\mathbb{R}\ $is a decreasing function,\ then we replace$\ \geq\ $by$\ \leq
\ $in (\ref{Ben}). \ 
\end{lemma}

\begin{remark}
The sign of $F$ plays no role in (\ref{Ben}).\ In the case when $M=S^{1},\ $it
is easy to obtain (\ref{Ben}) by Fubini theorem:
\begin{align*}
&  \int_{0}^{2\pi}d\theta\int_{0}^{2\pi}\xi\left(  \theta\right)  F\left(
\xi\left(  \theta\right)  \right)  d\theta-\int_{0}^{2\pi}\xi\left(
\theta\right)  d\theta\int_{0}^{2\pi}F\left(  \xi\left(  \theta\right)
\right)  d\theta\\
&  =\frac{1}{2}\int_{0}^{2\pi}\int_{0}^{2\pi}\left[  F\left(  \xi\left(
x\right)  \right)  -F\left(  \xi\left(  y\right)  \right)  \right]  \left[
\xi\left(  x\right)  -\xi\left(  y\right)  \right]  dxdy\geq0.
\end{align*}

\end{remark}

With the help of Andrews's inequality, we have\ (see p. 341 of \cite{A} also):

\begin{lemma}
(\emph{monotonicity of the IPD})\label{lem1}\ Under the general parabolic flow
(it could be contracting, expanding or a mixture of both)%
\begin{equation}
\dfrac{\partial X}{\partial t}\left(  \varphi,t\right)  =\left[  F\left(
k\left(  \varphi,t\right)  \right)  -\lambda\left(  t\right)  \right]
\mathbf{N}_{in}\left(  \varphi,t\right)  , \label{F5}%
\end{equation}
where $\lambda\left(  t\right)  \ $is a time function which may depend on the
global geometry $\gamma_{t}\ $of the flow.\ If\ the flow is well-defined on
$[0,T)$ and the evolving curves $\gamma_{t}$ stays embedded on $[0,T),\ $then
the IPD $L^{2}\left(  t\right)  -4\pi A\left(  t\right)  $ for $\gamma_{t}$ is
decreasing\ on $[0,T)$.
\end{lemma}

%

\proof
By (\ref{LA-evo}), we have
\begin{align}
&  \frac{d}{dt}\left(  L^{2}-4\pi A\right) \nonumber\\
&  =2L\left[  -\int_{\gamma_{t}}F\left(  k\right)  kds+2\pi\lambda\left(
t\right)  \right]  -4\pi\left[  -\int_{\gamma_{t}}F\left(  k\right)
ds+\lambda\left(  t\right)  L\right] \nonumber\\
&  =2\left[  \int_{\gamma_{t}}kds\int_{\gamma_{t}}F\left(  k\right)
ds-\int_{\gamma_{t}}ds\int_{\gamma_{t}}F\left(  k\right)  kds\right]  \leq0.
\label{09}%
\end{align}
In particular, we note that the function $\lambda\left(  t\right)  \ $has been
cancelled.\ The proof is done.$%
\hfill
\square$

\ \ \ 

As a consequence, we obtain the following:

\begin{corollary}
Under the assumption of Lemma \ref{lem1},\ if flow\ (\ref{F5}) is
area-preserving, then it must be length-decreasing.\ On the other hand,\ if it
is\ length-preserving, then it must be\ area-increasing.\ In particular, if it
preserves either area or length, then the IPR $L^{2}/4\pi A\ $is decreasing.\ 
\end{corollary}

\begin{remark}
We see that flows (\ref{F1})-(\ref{F3})\ are all IPR decreasing.\ 
\end{remark}

What happens to the IPR?\ We can compute%
\[
\frac{d}{dt}\left(  \frac{L^{2}}{4\pi A}\right)  =I+II,
\]
where%
\begin{align}
I  &  =\frac{2}{\left(  4\pi A\right)  ^{2}}\cdot\Phi,\nonumber\\
\Phi &  =\underline{L^{2}}\underbrace{\left(  \int_{\gamma_{t}}kds\right)
\left(  \int_{\gamma_{t}}F\left(  k\right)  ds\right)  }-\underline{4\pi
A}\underbrace{\left(  \int_{\gamma_{t}}ds\right)  \left(  \int_{\gamma_{t}%
}F\left(  k\right)  kds\right)  } \label{I}%
\end{align}
and
\begin{equation}
II=-\frac{4\pi\lambda\left(  t\right)  L}{\left(  4\pi A\right)  ^{2}}\left(
L^{2}-4\pi A\right)  =\left\{
\begin{array}
[c]{l}%
\leq0,\ \ \ \text{if\ }\lambda\left(  t\right)  \geq0%
\vspace{3mm}%
\\
\geq0,\ \ \ \text{if\ }\lambda\left(  t\right)  \leq0.
\end{array}
\right.  . \label{II}%
\end{equation}
In (\ref{I}), there is a competition between Andrews's inequality and the
isoperimetric inequality$\ L^{2}\geq4\pi A$. Hence it has no definite sign in general.\ 

By (\ref{II}) we also notice that, roughly speaking, the flow is
better-behaved if it tends to expand more\ ($\lambda\left(  t\right)  \geq0$),
and worse-behaved if it tends to contract more\ ($\lambda\left(  t\right)
\leq0$).

To explain the\emph{ gradient flow} of the IPR and IPD, by (\ref{W})\ and
(\ref{LA-evo}), we have
\begin{equation}
\frac{d}{dt}\left(  L^{2}-4\pi A\right)  \left(  t\right)  =-2L\int
_{\gamma_{t}}\left\langle W,\ \left(  k-\frac{2\pi}{L}\right)  \mathbf{N}%
_{in}\right\rangle ds \label{Q1}%
\end{equation}
and%
\begin{equation}
\frac{d}{dt}\left(  \frac{L^{2}}{4\pi A}\right)  \left(  t\right)  =-\frac
{L}{2\pi A}\int_{\gamma_{t}}\left\langle W,\ \left(  k-\frac{L}{2A}\right)
\mathbf{N}_{in}\right\rangle ds, \label{Q2}%
\end{equation}
where$\ W\left(  \varphi,t\right)  =\left(  \partial\gamma/\partial t\right)
\left(  \varphi,t\right)  $ is the speed vector of the flow.\ By
(\ref{Q1}),\ the nonlocal flow%
\begin{equation}
\dfrac{\partial X}{\partial t}\left(  \varphi,t\right)  =2L\left(  t\right)
\left[  k\left(  \varphi,t\right)  -\frac{2\pi}{L\left(  t\right)  }\right]
\mathbf{N}_{in}\left(  \varphi,t\right)  \label{F6}%
\end{equation}
is \emph{the gradient flow of the IPD functional.\ }It only differs from
Gage's flow (\ref{F1}) by a\emph{ }time factor $2L\left(  t\right)
.\ $If$\ X\left(  \varphi,t\right)  $ is a solution to Gage's flow (\ref{F1}),
the function%
\begin{equation}
\tilde{X}\left(  \varphi,\tau\right)  =X\left(  \varphi,t\left(  \tau\right)
\right)  \label{XX}%
\end{equation}
will then be a solution to the flow (\ref{F6}) if we choose $t\left(
\tau\right)  \ $to\ satisfy the identity%
\begin{equation}
\frac{dt}{d\tau}=2L\left(  t\right)  ,\ \ \ L\left(  t\right)
=\text{length\ of}\ X\left(  \varphi,t\right)  .
\end{equation}
To see this, by the chain rule%
\[
\dfrac{\partial\tilde{X}}{\partial\tau}\left(  \varphi,\tau\right)  =\frac
{dt}{d\tau}\dfrac{\partial X}{\partial t}\left(  \varphi,t\right)  =2L\left(
k-\frac{2\pi}{L}\right)  \mathbf{N}_{in}%
\]
and the relation%
\[
\tilde{L}\left(  \tau\right)  =L\left(  t\right)  ,\ \ \ \tilde{k}\left(
\varphi,\tau\right)  =k\left(  \varphi,t\right)  ,\ \ \ \mathbf{\tilde{N}%
}_{in}\left(  \varphi,\tau\right)  =\mathbf{N}_{in}\left(  \varphi,t\right)
,\ \ \ t=t\left(  \tau\right)
\]
we see that%
\begin{equation}
\dfrac{\partial\tilde{X}}{\partial\tau}\left(  \varphi,\tau\right)
=2\tilde{L}\left(  \tau\right)  \left[  \tilde{k}\left(  \varphi,\tau\right)
-\frac{2\pi}{\tilde{L}\left(  \tau\right)  }\right]  \mathbf{\tilde{N}}%
_{in}\left(  \varphi,\tau\right)  . \label{XY}%
\end{equation}
Thus we can say that the two flows (\ref{F1})\ and (\ref{F6}) are
\emph{equivalent}.\ 

Similarly, by (\ref{Q2}), the nonlocal flow%
\begin{equation}
\dfrac{\partial X}{\partial t}\left(  \varphi,t\right)  =\frac{L\left(
t\right)  }{2\pi A\left(  t\right)  }\left[  k\left(  \varphi,t\right)
-\frac{L\left(  t\right)  }{2A\left(  t\right)  }\right]  \mathbf{N}%
_{in}\left(  \varphi,t\right)  \label{F7}%
\end{equation}
is \emph{the gradient flow of the IPR functional.\ }Again, it only differs
from Jiang-Pan's flow (\ref{F2}) by a time factor $L\left(  t\right)  /2\pi
A\left(  t\right)  .\ $Thus the two flows (\ref{F2})\ and (\ref{F7}) are
\emph{equivalent}.

We summarize the following:

\begin{lemma}
(\emph{gradient flow of the IPR and IPD})\label{lem0}\ Gage's nonlocal flow
(\ref{F1}) is the gradient flow of the IPD functional,
and\ Jiang-Pan's\ nonlocal flow (\ref{F2}) is the gradient flow of the IPR
functional.\ Both flows decrease the IPR and IPD.
\end{lemma}

Regarding the isoperimetric\ behavior\ of a flow,\ another important
observation is\ the following:\ if we have $W=\mathbf{N}_{in}\ $in
(\ref{Q1})\ and\ (\ref{Q2}), then
\begin{equation}
\frac{d}{dt}\left(  L^{2}-4\pi A\right)  \left(  t\right)  =0,\ \ \ \ \ \frac
{d}{dt}\left(  \frac{L^{2}}{4\pi A}\right)  \left(  t\right)  =\frac{L}{4\pi
A^{2}}\left(  L^{2}-4\pi A\right)  \geq0 \label{Z1}%
\end{equation}
and if we have $W=u\mathbf{N}_{in},$ where $u=\left\langle \gamma\left(
\cdot,t\right)  ,\ \mathbf{N}_{out}\right\rangle $ is the \emph{support
function} of$\ \gamma\left(  \cdot,t\right)  ,$ then we have%
\begin{equation}
\frac{d}{dt}\left(  L^{2}-4\pi A\right)  \left(  t\right)  =-2\left(
L^{2}-4\pi A\right)  \leq0,\ \ \ \ \ \frac{d}{dt}\left(  \frac{L^{2}}{4\pi
A}\right)  \left(  t\right)  =0 \label{Z2}%
\end{equation}
due to the identities (see (\ref{kLA})\ also)%
\begin{equation}
\int_{\gamma_{t}}ukds=L,\ \ \ \ \ \frac{1}{2}\int_{\gamma_{t}}uds=A.
\end{equation}
Hence if one replace the flow
\begin{equation}
\dfrac{\partial X}{\partial t}\left(  \varphi,t\right)  =F\left(  k\left(
\varphi,t\right)  \right)  \mathbf{N}_{in}\left(  \varphi,t\right)
\label{XX0}%
\end{equation}
by$\ $%
\begin{equation}
\dfrac{\partial X}{\partial t}\left(  \varphi,t\right)  =\left[  F\left(
k\left(  \varphi,t\right)  \right)  -\lambda\left(  t\right)  \right]
\mathbf{N}_{in}\left(  \varphi,t\right)  \label{XX1}%
\end{equation}
then the IPD$\ L^{2}-4\pi A$ is unaffected. Similarly if one replaces
(\ref{XX0})\ by
\begin{equation}
\dfrac{\partial X}{\partial t}\left(  \varphi,t\right)  =\left[  F\left(
k\left(  \varphi,t\right)  \right)  -\lambda\left(  t\right)  u\left(
\varphi,t\right)  \right]  \mathbf{N}_{in}\left(  \varphi,t\right)  ,
\label{XX2}%
\end{equation}
where $u\left(  \varphi,t\right)  $ is the support function of $\gamma_{t}%
\ $at the point $\varphi,\ $then the IPR $L^{2}/4\pi A\ $is unaffected.\ 

Now we explain why the IPR $L^{2}/4\pi A\ $is unchanged by the extra term
$-\lambda u\mathbf{N}_{in}$.\ If $X\left(  \varphi,t\right)  $ is a solution
to (\ref{XX0}) and let $\tilde{X}\left(  \varphi,t\right)  =\sigma\left(
t\right)  X\left(  \varphi,t\right)  ,\ \sigma\left(  t\right)  >0,\ $which is
a time-dependent dilation of $X\left(  \varphi,t\right)  ,$ then the
dilated$\ \tilde{X}\left(  \varphi,t\right)  $ satisfies
\begin{equation}
\dfrac{\partial\tilde{X}}{\partial t}\left(  \varphi,t\right)  =\left\{
\begin{array}
[c]{l}%
\sigma\left(  t\right)  F\left(  \sigma\left(  t\right)  \tilde{k}\left(
\varphi,t\right)  \right)  \mathbf{\tilde{N}}_{in}\left(  \varphi,t\right)
\vspace{3mm}%
\\
+\dfrac{\sigma^{\prime}\left(  t\right)  }{\sigma\left(  t\right)
}\left\langle \tilde{X}\left(  \varphi,t\right)  ,\ \mathbf{\tilde{N}}%
_{out}\left(  \varphi,t\right)  \right\rangle \mathbf{\tilde{N}}_{out}\left(
\varphi,t\right)  +\dfrac{\sigma^{\prime}\left(  t\right)  }{\sigma\left(
t\right)  }\left\langle \tilde{X}\left(  \varphi,t\right)  ,\ \mathbf{\tilde
{T}}\left(  \varphi,t\right)  \right\rangle \mathbf{\tilde{T}}\left(
\varphi,t\right)  .
\end{array}
\right.
\end{equation}
Thus, up to a tangential component\ (it is known that a tangential component
can be removed by a further change of variable $\phi=\phi\left(
\varphi,t\right)  \ $in parametrizing the dilated curve $\tilde{\gamma}$), we
obtain
\begin{equation}
\dfrac{\partial\tilde{X}}{\partial t}\left(  \varphi,t\right)  =\left[  \sigma
F\left(  \sigma\tilde{k}\right)  -\dfrac{\sigma^{\prime}}{\sigma}\tilde
{u}\right]  \mathbf{\tilde{N}}_{in}, \label{XX3}%
\end{equation}
where the support function $\tilde{u}\ $appears naturally.\ Since a dilation
will not change the IPR,\ the extra term $-\sigma^{-1}\sigma^{\prime}\tilde
{u}$ must\ have no effect at all.\ 

In view of the above, one can\ keep dilating a solution$\ \gamma_{t}$ to the
flow (\ref{XX0}) so that its length or area is independent\ of time. The flow
equation for the dilated solution will have an extra term involving the
support function. There are two of them for the$\ k$-type flow$\ $%
with$\ F\left(  k\right)  =k,\ $and two of them for the$\ 1/k$-type
flow$\ $with$\ F\left(  k\right)  =-1/k.\ $In conclusion, we have the
following four flows:%
\begin{equation}
\dfrac{\partial X}{\partial t}\left(  \varphi,t\right)  =\left(  k-\frac{\pi
}{A}u\right)  \mathbf{N}_{in}\ \ \ \text{(area-preserving)} \label{S1}%
\end{equation}
and
\begin{equation}
\dfrac{\partial X}{\partial t}\left(  \varphi,t\right)  =\left[  k-\left(
\frac{1}{L}\int_{0}^{L}k^{2}ds\right)  u\right]  \mathbf{N}_{in}%
\ \ \ \text{(length-preserving)} \label{S2}%
\end{equation}
and%
\begin{equation}
\dfrac{\partial X}{\partial t}\left(  \varphi,t\right)  =\left[  \left(
\frac{1}{2A}\int_{0}^{L}\frac{1}{k}ds\right)  u-\frac{1}{k}\right]
\mathbf{N}_{in}\ \ \ \text{(area-preserving)} \label{S3}%
\end{equation}
and%
\begin{equation}
\dfrac{\partial X}{\partial t}\left(  \varphi,t\right)  =\left(  u-\frac{1}%
{k}\right)  \mathbf{N}_{in}\ \ \ \text{(length-preserving),} \label{S4}%
\end{equation}
where $u\left(  \varphi,t\right)  $ is the support function of the curve
$\gamma_{t}.\ $For example, for $F\left(  k\right)  =k,\ $(\ref{XX3}) becomes%
\[
\dfrac{\partial X}{\partial t}\left(  \varphi,t\right)  =\sigma^{2}\left(
k-\dfrac{\sigma^{\prime}}{\sigma^{3}}u\right)  \mathbf{N}_{in}%
\]
and a change in time variable can get rid of the coefficient$\ \sigma^{2}%
.\ $Hence we may assume
\begin{equation}
\dfrac{\partial X}{\partial t}\left(  \varphi,t\right)  =\left(
k-\dfrac{\sigma^{\prime}}{\sigma^{3}}u\right)  \mathbf{N}_{in}. \label{S5}%
\end{equation}
If we want to keep the enclosed area fixed, by (\ref{LA-evo})\ we need to
require$\ $%
\[
\frac{dA}{dt}\left(  t\right)  =-\int_{\gamma_{t}}\left(  k-\dfrac
{\sigma^{\prime}}{\sigma^{3}}u\right)  ds=0,
\]
which implies$\ \sigma^{\prime}\left(  t\right)  /\sigma^{3}\left(  t\right)
=\pi/A\left(  t\right)  $ and the flow (\ref{S5})\ becomes
\begin{equation}
\dfrac{\partial X}{\partial t}\left(  \varphi,t\right)  =\left(  k-\dfrac{\pi
}{A}u\right)  \mathbf{N}_{in}.
\end{equation}
The same argument can be applied to the other three flows.\ 

We find that if we replace $u$ by $2A/L\ $in the area-preserving
flows\ (\ref{S1}) and (\ref{S3}),\ we get Gage's flow (\ref{F1})\ and
Ma-Cheng's flow (\ref{F4}).\ Also if we replace$\ u$ by $L/2\pi\ $in the
length-preserving flows\ (\ref{S2}) and (\ref{S4}),\ we get\ Ma-Zhu's flow
(\ref{F2-1}) and\ Pan-Yang's flow (\ref{F3}).\ This is due to formula
(\ref{kLA}) below.\ 

\section{The\ dual flow of (\ref{F2}) and an improved isoperimetric
inequality\label{sec}}

We already\ know that flows\ (\ref{F1})\ and (\ref{F3})\ are dual to each
other. Motivated by it, we can also consider the dual flow of (\ref{F2}),
which is probably the only remaining interesting case not dealt with among
those nonlocal flows (\ref{F1})-(\ref{F3}).\ It has the form%
\begin{equation}
\left\{
\begin{array}
[c]{l}%
\dfrac{\partial X}{\partial t}\left(  \varphi,t\right)  =\left[
\dfrac{2A\left(  t\right)  }{L\left(  t\right)  }-\dfrac{1}{k\left(
\varphi,t\right)  }\right]  \mathbf{N}_{in}\left(  \varphi,t\right)
\vspace{3mm}%
\\
X\left(  \varphi,0\right)  =X_{0}\left(  \varphi\right)  ,\ \ \ \varphi\in
S^{1},
\end{array}
\right.  . \label{dual1}%
\end{equation}
where$\ X_{0}\left(  \varphi\right)  $ is the parametrization of a given
smooth\ convex closed curve $\gamma_{0}.\ $

In $k$-type\ flows, the Gage's inequality for a\ convex closed\ curve
$\gamma:$
\begin{equation}
\int_{\gamma}k^{2}\left(  s\right)  ds=\int_{0}^{2\pi}k\left(  \theta\right)
d\theta\geq\frac{\pi L}{A}, \label{Gage}%
\end{equation}
where $s$ is arc\ length\ parameter, plays an important role.\ In contrast, in
the $1/k$-type\ flows, we need to use the Pan-Yang's isoperimetric inequality
for a\ convex closed\ curve $\gamma:$
\begin{equation}
\int_{\gamma}\frac{1}{k\left(  s\right)  }ds\geq\frac{L^{2}-2\pi A}{\pi},
\label{Pan}%
\end{equation}
where the equality holds if and only if $\gamma$ is a circle.\ 

(\ref{Pan}) is proved in \cite{PY} using methods\ established in Green-Osher
\cite{GO}.\ This inequality seems to be new as it had not appeared in any book
or reference before.\ Here we can use Fourier series expansion\ to give an
alternative proof and,\ at the same time, improves it also.\ 

In the book by Courant-John \cite{CJ}, they used support function and Fourier
series method to prove the isoperimetric inequality $L^{2}\geq4\pi A\ $for a
closed plane curve (see p. 366 of \cite{CJ}).\ Our method is motivated by
theirs.\ Let $C$ be a convex closed plane\ curve.\ One can use its outward
normal angle $\theta\in\left[  0,2\pi\right]  $ to parametrize it.\ In doing
so, the inequality (\ref{Pan})\ becomes%
\begin{equation}
\int_{0}^{2\pi}\frac{1}{k^{2}\left(  \theta\right)  }d\theta\geq\frac
{L^{2}-2\pi A}{\pi}. \label{Pan1}%
\end{equation}
It is also known that one can use the support function $u\left(
\theta\right)  ,\ \theta\in\left[  0,2\pi\right]  ,\ $of $C$ to express its
curvature, enclosed area and length\ (see the book by\ Schneider \cite{S}). We
have%
\begin{equation}
\left\{
\begin{array}
[c]{l}%
\dfrac{1}{k\left(  \theta\right)  }=u_{\theta\theta}\left(  \theta\right)
+u\left(  \theta\right)  ,\ \ \ L=%
{\displaystyle\int_{0}^{2\pi}}
u\left(  \theta\right)  d\theta%
\vspace{3mm}%
\\
A=\dfrac{1}{2}%
{\displaystyle\int_{C}}
uds=\dfrac{1}{2}%
{\displaystyle\int_{0}^{2\pi}}
u\left(  \theta\right)  \left[  u_{\theta\theta}\left(  \theta\right)
+u\left(  \theta\right)  \right]  d\theta.
\end{array}
\right.  \label{kLA}%
\end{equation}
We can state our result in the following:

\begin{lemma}
(\emph{refined Pan-Yang's isoperimetric\ inequality})\ For any$\ $convex
closed plane\ curve\ $C$ there holds the inequality%
\begin{equation}
\int_{C}\frac{1}{k\left(  s\right)  }ds=\int_{0}^{2\pi}\frac{1}{k^{2}\left(
\theta\right)  }d\theta\geq\frac{2}{\pi}\left(  L^{2}-4\pi A\right)  +2A,
\label{Tsai}%
\end{equation}
where the equality holds if and only if the support function of $C\ $has the
form%
\begin{equation}
u\left(  \theta\right)  =a_{0}+a_{1}\cos\theta+b_{1}\sin\theta+a_{2}%
\cos2\theta+b_{2}\sin2\theta,\ \ \ \theta\in\left[  0,2\pi\right]
\label{const}%
\end{equation}
for some constants $a_{0},\ a_{1},\ b_{1},\ a_{2},\ b_{2}\ $satisfying%
\begin{equation}
u_{\theta\theta}\left(  \theta\right)  +u\left(  \theta\right)  =a_{0}%
-3a_{2}\cos2\theta-3b_{2}\sin2\theta>0\ \ \ \text{for all\ \ \ }\theta
\in\left[  0,2\pi\right]  .
\end{equation}
Here the variable $\theta$ is the outward normal angle of $C.\ $
\end{lemma}

\begin{remark}
Since
\[
\frac{2}{\pi}\left(  L^{2}-4\pi A\right)  +2A=\frac{L^{2}-4\pi A}{\pi}%
+\frac{L^{2}-2\pi A}{\pi}\geq\frac{L^{2}-2\pi A}{\pi},
\]
(\ref{Tsai}) is an improvement of (\ref{Pan}).\ 
\end{remark}

%

\proof
Using Fourier series, one can express the $2\pi$-periodic smooth support
function$\ u\left(  \theta\right)  $ of $C\ $as%
\begin{equation}
u\left(  \theta\right)  =\frac{a_{0}}{2}+\sum_{n=1}^{\infty}\left(  a_{n}\cos
n\theta+b_{n}\sin n\theta\right)  . \label{08-2}%
\end{equation}
Then%
\begin{equation}
L^{2}=\left(  \int_{0}^{2\pi}ud\theta\right)  ^{2}=\left(  \pi a_{0}\right)
^{2} \label{08-1}%
\end{equation}
and
\begin{equation}
4\pi A=\left(  \pi a_{0}\right)  ^{2}+2\pi^{2}\left[  \sum_{n=1}^{\infty
}\left(  1-n^{2}\right)  \left(  a_{n}^{2}+b_{n}^{2}\right)  \right]
=L^{2}+2\pi^{2}\left[  \sum_{n=2}^{\infty}\left(  1-n^{2}\right)  \left(
a_{n}^{2}+b_{n}^{2}\right)  \right]  , \label{08}%
\end{equation}
which gives the classical isoperimetric inequality$\ L^{2}\geq4\pi A.$\ Also
we have%
\[
\int_{0}^{2\pi}u_{\theta\theta}\left(  u_{\theta\theta}+u\right)  d\theta
=\sum_{n=2}^{\infty}n^{2}\left(  n^{2}-1\right)  \pi\left(  a_{n}^{2}%
+b_{n}^{2}\right)
\]
and therefore%
\begin{align*}
&  \int_{0}^{2\pi}\frac{1}{k^{2}\left(  \theta\right)  }d\theta\\
&  =\int_{0}^{2\pi}u_{\theta\theta}\left(  u_{\theta\theta}+u\right)
d\theta+\int_{0}^{2\pi}u\left(  u_{\theta\theta}+u\right)  d\theta=\sum
_{n=2}^{\infty}n^{2}\left(  n^{2}-1\right)  \pi\left(  a_{n}^{2}+b_{n}%
^{2}\right)  +2A.
\end{align*}
To prove (\ref{Tsai}), it suffices to show that%
\[
\sum_{n=2}^{\infty}n^{2}\left(  n^{2}-1\right)  \pi\left(  a_{n}^{2}+b_{n}%
^{2}\right)  \geq\frac{2}{\pi}\left(  L^{2}-4\pi A\right)  ,
\]
where by (\ref{08})\ the right hand side is%
\[
\frac{2}{\pi}\left(  L^{2}-4\pi A\right)  =4\pi\left[  \sum_{n=2}^{\infty
}\left(  n^{2}-1\right)  \left(  a_{n}^{2}+b_{n}^{2}\right)  \right]  .
\]
We clearly have%
\begin{equation}
\sum_{n=2}^{\infty}n^{2}\left(  n^{2}-1\right)  \pi\left(  a_{n}^{2}+b_{n}%
^{2}\right)  \geq4\pi\left[  \sum_{n=2}^{\infty}\left(  n^{2}-1\right)
\left(  a_{n}^{2}+b_{n}^{2}\right)  \right]
\end{equation}
and the equality holds if and only if\ $a_{n}=b_{n}=0\ $for all $n\geq
3.\ $That is, if and only if $u\left(  \theta\right)  $ has the form%
\[
u\left(  \theta\right)  =c+a_{1}\cos\theta+b_{1}\sin\theta+a_{2}\cos
2\theta+b_{2}\sin2\theta,\ \ \ \theta\in\left[  0,2\pi\right]
\]
for some constants $c,\ a_{1},\ b_{1},\ a_{2},\ b_{2}.\ $The proof is done.$%
\hfill
\square$

\begin{remark}
It is very unlikely to use Fourier series expansion method to prove Gage's
inequality (\ref{Gage})\ because the integrand$\ k\left(  \theta\right)
=\left[  u_{\theta\theta}\left(  \theta\right)  +u\left(  \theta\right)
\right]  ^{-1}\ $is not of the right form.\ 
\end{remark}

Regarding flows (\ref{F1})-(\ref{F3}) and (\ref{dual1}), I think one can use
the result in section 2 of Gage-Hamilton \cite{GH} (\emph{Nash-Moser inverse
function} \emph{theorem})\ to prove that for any initial curve, the nonlocal
flow considered has a solution for short time.\ However, since the\ Nash-Moser
inverse function theorem is itself hard to understand, one would prefer to
invoke a more straightforward or elementary\ theory.\ As the dual
flow\ (\ref{dual1}) has more of a \emph{linear structure}\ (for the support
function or for the inverse of curvature),\ it is possible to prove short time
existence directly using Fourier series.\ 

In Pan-Yang's length-preserving$\ $flow (\ref{F3}), they pointed out that if
the flow has a smooth convex\ solution on short time interval$\ [0,T),\ $then
the function (in below the variable$\ \theta$ represents outward normal angle
of the convex curve$\ \gamma_{t}$)%
\begin{equation}
w\left(  \theta,t\right)  :=e^{-t}\left(  \frac{1}{k\left(  \theta,t\right)
}-\frac{L\left(  t\right)  }{2\pi}\right)  ,\ \ \ L\left(  t\right)  =L\left(
0\right)  ,\ \ \ \left(  \theta,t\right)  \in S^{1}\times\lbrack0,T)
\label{07}%
\end{equation}
will satisfy a standard linear heat equation $w_{t}\left(  \theta,t\right)
=w_{\theta\theta}\left(  \theta,t\right)  \ $on $S^{1}\ $with initial
condition $w_{0}\left(  \theta\right)  =1/k_{0}\left(  \theta\right)
-L\left(  0\right)  /2\pi,\ k_{0}\left(  \theta\right)  >0.\ $Conversely one
can use the linear\ heat equation to establish short time existence\ of a
solution to the nonlocal flow\ (\ref{F2}) because if one knows$\ w\left(
\theta,t\right)  $ then one\ can know$\ $the curvature$\ k\left(
\theta,t\right)  \ $(here we need to use the fact $L\left(  t\right)
=L\left(  0\right)  \ $is preserved), and then we use curvature to construct a
solution to the nonlocal flow\ (\ref{F3}).

For the dual flow (\ref{dual1}),\ the situation is different. If it has a
smooth\ convex solution on $[0,T)\ $for short time$\ T>0,$ then the curvature
and length satisfy the following evolution equations:\ \
\begin{equation}
\frac{\partial}{\partial t}\left(  \frac{1}{k\left(  \theta,t\right)
}\right)  =\left(  \frac{1}{k\left(  \theta,t\right)  }\right)  _{\theta
\theta}+\frac{1}{k\left(  \theta,t\right)  }-\dfrac{2A\left(  t\right)
}{L\left(  t\right)  } \label{06}%
\end{equation}
and
\begin{equation}
\frac{d}{dt}\left(  \frac{L\left(  t\right)  }{2\pi}\right)  =-\frac{1}{2\pi}%
{\displaystyle\int_{\gamma_{t}}}
\left(  \dfrac{2A\left(  t\right)  }{L\left(  t\right)  }-\dfrac{1}{k}\right)
kds=\frac{L\left(  t\right)  }{2\pi}-\dfrac{2A\left(  t\right)  }{L\left(
t\right)  }\geq0. \label{05}%
\end{equation}
Hence we have the nice-looking equation%
\begin{equation}
\frac{\partial}{\partial t}\left(  \frac{1}{k\left(  \theta,t\right)  }%
-\frac{L\left(  t\right)  }{2\pi}\right)  =\left(  \frac{1}{k\left(
\theta,t\right)  }-\frac{L\left(  t\right)  }{2\pi}\right)  _{\theta\theta
}+\left(  \frac{1}{k\left(  \theta,t\right)  }-\frac{L\left(  t\right)  }%
{2\pi}\right)  \label{01}%
\end{equation}
and the function
\begin{equation}
w\left(  \theta,t\right)  :=e^{-t}\left(  \frac{1}{k\left(  \theta,t\right)
}-\frac{L\left(  t\right)  }{2\pi}\right)  ,\ \ \ \frac{L\left(  t\right)
}{2\pi}=\frac{1}{2\pi}\int_{0}^{2\pi}\frac{1}{k\left(  \theta,t\right)
}d\theta
\end{equation}
satisfies a \emph{linear\ heat equation}$\ $%
\begin{equation}
w_{t}\left(  \theta,t\right)  =w_{\theta\theta}\left(  \theta,t\right)
,\ \ \ w\left(  \theta,0\right)  =\frac{1}{k_{0}\left(  \theta\right)  }%
-\frac{L\left(  0\right)  }{2\pi}. \label{heat-eq}%
\end{equation}
Unfortunately, if we go in the reverse direction and solve $w\left(
\theta,t\right)  $ from the heat equation\ (\ref{heat-eq}), we are \emph{not}
able to recover the curvature $k\left(  \theta,t\right)  $ since from an
identity of the form
\begin{equation}
\frac{1}{k\left(  \theta,t\right)  }-\frac{1}{2\pi}\int_{0}^{2\pi}\frac
{1}{k\left(  \theta,t\right)  }d\theta=e^{t}w\left(  \theta,t\right)
\label{fw}%
\end{equation}
we cannot determine $k\left(  \theta,t\right)  \ $uniquely$\ $(note that\ if
$1/k\left(  \theta,t\right)  \ $satisfies (\ref{fw}),\ so is $1/k\left(
\theta,t\right)  +g\left(  t\right)  \ $for any\ function of time $g\left(
t\right)  $). In addition, even we find one $k\left(  \theta,t\right)
\ $satisfies (\ref{fw}),\ we do not know if it satisfies the evolution
equation (\ref{06}).\ The same difficulty also happens in Ma-Cheng's flow
(\ref{F4}). See Theorem 16 of \cite{MC}.

\begin{remark}
Note that in Pan-Yang's flow (\ref{F3}) we have $L\left(  t\right)  =L\left(
0\right)  $ for all $t.\ $Hence if $w\left(  \theta,t\right)  $ is known
from\ (\ref{heat-eq}), one can determine $k\left(  \theta,t\right)  \ $and its
evolution equation, and the short time existence of a solution to the flow is established.\ 
\end{remark}

To overcome the above difficulty, we can use Fourier series method again.\ The
evolution of the support function $u\left(  \theta,t\right)  \ $under flow
(\ref{dual1})\ on $[0,T)\ $is given by
\begin{equation}
\frac{\partial u}{\partial t}\left(  \theta,t\right)  =u_{\theta\theta}\left(
\theta,t\right)  +u\left(  \theta,t\right)  -\dfrac{\int_{0}^{2\pi}u\left(
\theta,t\right)  \left(  u_{\theta\theta}\left(  \theta,t\right)  +u\left(
\theta,t\right)  \right)  d\theta}{\int_{0}^{2\pi}u\left(  \theta,t\right)
d\theta},\ \ \ \left(  \theta,t\right)  \in S^{1}\times\lbrack0,T) \label{dtu}%
\end{equation}
where $u\left(  \theta,0\right)  =u_{0}\left(  \theta\right)  ,\ u_{0}%
^{\prime\prime}\left(  \theta\right)  +u_{0}\left(  \theta\right)
>0,\ $and$\ u_{0}\left(  \theta\right)  \ $is the support function of the
initial convex curve $\gamma_{0}.\ $One can use relation (\ref{kLA}) to derive
the equation (\ref{dtu}).\ We want to use (\ref{dtu}), instead of equation
(\ref{heat-eq}), to prove the existence of a nonlocal\ flow solution (the fact
is that if we have a support function\ solution to (\ref{dtu}), one can use it
to construct a flow solution to (\ref{dual1})).

Expand $u_{0}\left(  \theta\right)  \ $and$\ u\left(  \theta,t\right)  $ as%
\[
\left\{
\begin{array}
[c]{l}%
u_{0}\left(  \theta\right)  =\dfrac{a_{0}\left(  0\right)  }{2}+%
{\displaystyle\sum_{n=1}^{\infty}}
\left(  a_{n}\left(  0\right)  \cos n\theta+b_{n}\left(  0\right)  \sin
n\theta\right)
\vspace{3mm}%
\\
u\left(  \theta,t\right)  =\dfrac{a_{0}\left(  t\right)  }{2}+%
{\displaystyle\sum_{n=1}^{\infty}}
\left(  a_{n}\left(  t\right)  \cos n\theta+b_{n}\left(  t\right)  \sin
n\theta\right)  ,\ \ \ u\left(  \theta,0\right)  =u_{0}\left(  \theta\right)
\end{array}
\right.
\]
and compute%
\begin{align}
\frac{\partial}{\partial t}u\left(  \theta,t\right)   &  =\frac{a_{0}^{\prime
}\left(  t\right)  }{2}+\sum_{n=1}^{\infty}\left(  a_{n}^{\prime}\left(
t\right)  \cos n\theta+b_{n}^{\prime}\left(  t\right)  \sin n\theta\right)
\nonumber\\
u_{\theta\theta}\left(  \theta,t\right)  +u\left(  \theta,t\right)   &
=\frac{a_{0}\left(  t\right)  }{2}+\sum_{n=1}^{\infty}\left(  1-n^{2}\right)
\left(  a_{n}\left(  t\right)  \cos n\theta+b_{n}\left(  t\right)  \sin
n\theta\right)  \label{star}%
\end{align}
and%
\[
\dfrac{\int_{0}^{2\pi}u\left(  \theta,t\right)  \left(  u_{\theta\theta
}\left(  \theta,t\right)  +u\left(  \theta,t\right)  \right)  d\theta}%
{\int_{0}^{2\pi}u\left(  \theta,t\right)  d\theta}=\frac{\frac{\pi}{2}%
a_{0}^{2}\left(  t\right)  +\pi\sum_{n=1}^{\infty}\left(  1-n^{2}\right)
\left(  a_{n}^{2}\left(  t\right)  +b_{n}^{2}\left(  t\right)  \right)  }{\pi
a_{0}\left(  t\right)  }.
\]
Hence, by comparing the coefficients, we want $a_{0}\left(  t\right)
,\ a_{n}\left(  t\right)  ,\ b_{n}\left(  t\right)  $ to satisfy
\begin{equation}
\left\{
\begin{array}
[c]{l}%
\dfrac{a_{0}^{\prime}\left(  t\right)  }{2}=\dfrac{\sum_{n=1}^{\infty}\left(
n^{2}-1\right)  \left(  a_{n}^{2}\left(  t\right)  +b_{n}^{2}\left(  t\right)
\right)  }{a_{0}\left(  t\right)  }\geq0,\ \ \ a_{0}\left(  0\right)  >0%
\vspace{3mm}%
\\
a_{n}^{\prime}\left(  t\right)  =\left(  1-n^{2}\right)  a_{n}\left(
t\right)  ,\ \ \ \ \ b_{n}^{\prime}\left(  t\right)  =\left(  1-n^{2}\right)
b_{n}\left(  t\right)  .
\end{array}
\right.  \label{ODE}%
\end{equation}
This ODE system can be easily solved for all $t\in\lbrack0,\infty)\ $to get%
\begin{equation}
\left\{
\begin{array}
[c]{l}%
a_{n}\left(  t\right)  =a_{n}\left(  0\right)  e^{\left(  1-n^{2}\right)
t},\ \ \ \ \ b_{n}\left(  t\right)  =b_{n}\left(  0\right)  e^{\left(
1-n^{2}\right)  t},\ \ \ \ \ n\in\mathbb{N}%
\vspace{3mm}%
\\
a_{0}\left(  t\right)  =\sqrt{a_{0}^{2}\left(  0\right)  +2%
{\displaystyle\sum_{n=1}^{\infty}}
\left(  1-e^{2\left(  1-n^{2}\right)  t}\right)  \left(  a_{n}^{2}\left(
0\right)  +b_{n}^{2}\left(  0\right)  \right)  }.
\end{array}
\right.  \label{ODE-sol}%
\end{equation}
Finally we note that as $a_{0}\left(  t\right)  ,\ a_{n}\left(  t\right)
,\ b_{n}\left(  t\right)  $ are all exponentially decay,\ the above
computations can all be justified by Fourier series theory.

In view of the above, we can conclude the following:

\begin{lemma}
(\emph{short time existence of the dual flow (\ref{dual1})})\ There is a
smooth\ convex solution $X\left(  \varphi,t\right)  $ to the nonlocal flow
(\ref{dual1}) for short time interval$\ [0,T),\ T>0.\ $
\end{lemma}

A major difficulty in studying the $1/k$-type\ flows (\ref{F3}),\ (\ref{F4}),
(\ref{dual1})$\ $is the possibility of developing a singularity in
\emph{finite} time with curvature $k=\infty$ somewhere.\ Although this seems
quite unlikely to happen, we are not able to rule it out mathematically.\ Note
that the flow (\ref{dual1}) is equivalent\ to the support
function\ equation\ (\ref{dtu}) only under the condition
\begin{equation}
0<u_{\theta\theta}\left(  \theta,t\right)  +u\left(  \theta,t\right)
<\infty\label{uu}%
\end{equation}
since the curvature is given by $k=1/\left(  u_{\theta\theta}+u\right)  .$
From equation\ (\ref{dtu}) we see that $u_{\theta\theta}\left(  \theta
,t\right)  +u\left(  \theta,t\right)  $ will not blow up in finite time, but
it may be possible that $u_{\theta\theta}\left(  \theta_{0},t_{0}\right)
+u\left(  \theta_{0},t_{0}\right)  =0\ $at some finite time $t_{0}$ for some
$\theta_{0}\in S^{1}\ $(note that the initial condition satisfies$\ u_{0}%
^{\prime\prime}\left(  \theta\right)  +u_{0}\left(  \theta\right)  >0$
everywhere). Even we have an explicit Fourier series expansion for
$u_{\theta\theta}\left(  \theta,t\right)  +u\left(  \theta,t\right)  ,$ we do
not know how to exclude the possibility.\ 

\begin{remark}
In both \cite{PY} and \cite{MC}, although they claim that the flow they
studied will converge to a round circle as $t\rightarrow\infty,\ $the
possibility of $k\ $becoming infinity in finite time is not
discussed.\ However, this should \emph{not} diminish their contributions to
the study of nonlocal flows. In both papers, they derived the estimate (see
\cite{PY}, p. 481\ and \cite{MC}, p. 8)%
\begin{equation}
\left\vert \frac{1}{k\left(  \theta,t\right)  }-\frac{L\left(  t\right)
}{2\pi}\right\vert \leq Me^{T^{\ast}}\ \ \ \text{for any\ finite\ }T^{\ast}>0.
\end{equation}
This can exclude a finite time extinction ($k=0$)\ of the curvature, but it
can \textbf{not} exclude a finite time blow-up ($k=\infty$)\ of the curvature.\ 
\end{remark}

For the dual flow (\ref{dual1}),\ if we only look at equation (\ref{dtu}), we
can easily obtain the following convergence. As mentioned above, it can give
information of the flow only when (\ref{uu}) is satisfied.\ 

\begin{theorem}
\label{thm1}Let$\ u_{0}\left(  \theta\right)  $ be a smooth function on
$S^{1}$ satisfying $u_{0}^{\prime\prime}\left(  \theta\right)  +u_{0}\left(
\theta\right)  >0$ everywhere.\ Then the solution to the equation (\ref{dtu})
with initial condition $u\left(  \theta,0\right)  =u_{0}\left(  \theta\right)
,\ \theta\in S^{1},$ is defined on $S^{1}\times\lbrack0,\infty)$ with
\begin{equation}
\lim_{t\rightarrow\infty}\left\Vert u\left(  \theta,t\right)  -\left(
c+a_{1}\left(  0\right)  \cos\theta+b_{1}\left(  0\right)  \sin\theta\right)
\right\Vert _{C^{k}\left(  S^{1}\right)  }=0\ \ \ \text{for any\ \ \ }%
k\in\mathbb{N}%
\end{equation}
where $c>0$ is a constant given by%
\begin{align}
c  &  =\frac{1}{2}\sqrt{a_{0}^{2}\left(  0\right)  +2%
{\displaystyle\sum_{n=2}^{\infty}}
\left(  a_{n}^{2}\left(  0\right)  +b_{n}^{2}\left(  0\right)  \right)
}\nonumber\\
&  =\frac{1}{\sqrt{2\pi}}\left(  \int_{0}^{2\pi}\left[  u_{0}\left(
\theta\right)  -a_{1}\left(  0\right)  \cos\theta-b_{1}\left(  0\right)
\sin\theta\right]  ^{2}d\theta\right)  ^{1/2}. \label{c}%
\end{align}
Here $a_{0}\left(  0\right)  ,\ a_{n}\left(  0\right)  ,\ b_{n}\left(
0\right)  $ are the Fourier coefficients of the function $u_{0}\left(
\theta\right)  .\ $
\end{theorem}

%

\proof
This is obvious from Fourier series expansion.\ Perhaps one would worry about
the denominator in (\ref{dtu}) being zero in finite time, but this will not
happen from the Fourier series expansion of $a_{0}\left(  t\right)  ,\ $or one
can compute$\ $(note that$\ \int_{0}^{2\pi}u_{0}\left(  \theta\right)
d\theta>0$)$\ $%
\[
\frac{d}{dt}%
{\displaystyle\int_{0}^{2\pi}}
u\left(  \theta,t\right)  d\theta=\int_{0}^{2\pi}u\left(  \theta,t\right)
d\theta-\dfrac{2\pi\int_{0}^{2\pi}u\left(  \theta,t\right)  \left(
u_{\theta\theta}\left(  \theta,t\right)  +u\left(  \theta,t\right)  \right)
d\theta}{\int_{0}^{2\pi}u\left(  \theta,t\right)  d\theta}\geq0
\]
due to the \emph{Poincar\'{e} inequality }for$\ 2\pi$-periodic
functions\footnote{It is interesting to know that, in fact, the Poincar\'{e}
inequality,\ the classical isoperimetric inequality$\ L^{2}\geq4\pi A$, and
the\ Minkowski\ mixed area\ inequality$\ \sqrt{A_{1}A_{2}}\leq A_{12}$\ (see
(\ref{AA})\ below) are all equivalent.\ }\ (see p. 179 of the book
\cite{BGH})\ :%
\begin{equation}
\int_{0}^{2\pi}dx\int_{0}^{2\pi}f\left(  \frac{d^{2}f}{dx^{2}}+f\right)
dx\leq\left(  \int_{0}^{2\pi}f\left(  x\right)  dx\right)  ^{2}. \label{Wir}%
\end{equation}
$%
\hfill
\square$

\ \ \ 

We now assume that the dual flow (\ref{dual1}) will not develop a singularity
($k=\infty$)\ in finite time.\ Then the flow must be defined on the infinite
time interval $[0,\infty),\ $with$\ $each\ $\gamma_{t}$ remaining smooth and
convex,\ and then we can look at its asymptotic\ geometry.\ The convexity of
$\gamma_{t}$ can be seen from equation (\ref{06}) since $2A\left(  t\right)
/L\left(  t\right)  \ $is uniformly bounded (see below)\ and by the maximum
principle, $1/k\left(  \theta,t\right)  $ will not blow up in finite time,
which implies that $k\left(  \theta,t\right)  >0$ will not become zero in
finite time.

In below\ we can quickly prove the convergence of the flow again without
relying on the support function and its Fourier series expansion.\ The
evolution of the length $L\left(  t\right)  $ is known by (\ref{05}).\ As for
area $A\left(  t\right)  $, we have
\begin{equation}
\frac{dA}{dt}=-2A+%
{\displaystyle\int_{\gamma_{t}}}
\frac{1}{k}ds\geq-2A+\frac{L^{2}-2\pi A}{\pi}=\frac{L^{2}-4\pi A}{\pi}\geq0,
\label{04}%
\end{equation}
where we have used the Pan-Yang's isoperimetric inequality (\ref{Pan}) in
(\ref{04}).\ Hence in flow (\ref{dual1}), \emph{both length and area are
increasing}. In particular%
\[
\frac{d}{dt}\left(  \dfrac{A}{L}\right)  \geq\dfrac{\left(  L^{2}-4\pi
A\right)  \left(  L^{2}-\pi A\right)  }{\pi L^{3}}\geq0.
\]
As a consequence, we obtain
\begin{equation}
\frac{d}{dt}\left(  L^{2}-4\pi A\right)  \leq2L\frac{L^{2}-4\pi A}{L}%
-4\pi\frac{L^{2}-4\pi A}{\pi}=-2\left(  L^{2}-4\pi A\right)  \leq0,
\label{04-1}%
\end{equation}
and derive the\emph{ }exponential decay\emph{ }of the IPD%
\begin{equation}
0\leq L^{2}\left(  t\right)  -4\pi A\left(  t\right)  \leq e^{-2t}\left(
L^{2}\left(  0\right)  -4\pi A\left(  0\right)  \right)  \rightarrow
0\ \ \ \text{as\ \ \ }t\rightarrow\infty.
\end{equation}
In particular, the IPR%
\begin{equation}
1\leq\frac{L^{2}}{4\pi A}=\frac{L^{2}-4\pi A}{4\pi A}+1\leq\frac
{e^{-2t}\left(  L^{2}\left(  0\right)  -4\pi A\left(  0\right)  \right)
}{4\pi A\left(  0\right)  }+1\rightarrow1\ \ \ \text{as\ \ \ }t\rightarrow
\infty.
\end{equation}
is exponentially decaying\textbf{ }to $1\ $since $A\left(  t\right)  \ $is increasing.

Due to the exponential decay of $L^{2}\left(  t\right)  -4\pi A\left(
t\right)  ,\ $the increasing $L\left(  t\right)  \ $and $A\left(  t\right)
\ $actually converge\ as $t\rightarrow\infty$.\ By%
\begin{equation}
0\leq\frac{d}{dt}\left(  \frac{L^{2}}{2}\right)  =L\frac{dL}{dt}=L^{2}-4\pi
A\leq e^{-2t}\left(  L^{2}\left(  0\right)  -4\pi A\left(  0\right)  \right)
\end{equation}
we have
\[
0\leq\frac{L^{2}\left(  t\right)  }{2}-\frac{L^{2}\left(  0\right)  }{2}%
\leq\frac{1}{2}\left[  L^{2}\left(  0\right)  -4\pi A\left(  0\right)
\right]  \left(  1-e^{-2t}\right)  .
\]
In particular, the limit$\ \lim_{t\rightarrow\infty}L\left(  t\right)
=L\left(  \infty\right)  >0\ $exists\ (note that $L\left(  t\right)  \ $is
increasing), where
\begin{equation}
L\left(  0\right)  \leq L\left(  \infty\right)  \leq\sqrt{L^{2}\left(
0\right)  +\left[  L^{2}\left(  0\right)  -4\pi A\left(  0\right)  \right]  }.
\end{equation}
As for $A\left(  t\right)  ,\ $we have$\ \lim_{t\rightarrow\infty}A\left(
t\right)  =A\left(  \infty\right)  =L^{2}\left(  \infty\right)  /4\pi.$ We
also have $\lim_{t\rightarrow\infty}\left(  2A\left(  t\right)  /L\left(
t\right)  \right)  =2A\left(  \infty\right)  /L\left(  \infty\right)  $ in
(\ref{06}).$\ $

\begin{remark}
Note that if we do not use Pan-Yang's isoperimetric inequality (\ref{Pan}) in
(\ref{04-1}), by H\"{o}lder inequality\ we would only obtain%
\[
\frac{d}{dt}\left(  L^{2}-4\pi A\right)  =2\left(  L^{2}-2\pi%
{\displaystyle\int_{\gamma_{t}}}
\frac{1}{k}ds\right)  \leq0,
\]
which is not good enough to imply the convergence of $L\left(  t\right)  \ $as
$t\rightarrow\infty.$
\end{remark}

To show the convergence of $1/k,\ $we can apply the following
simple\ result\ (see \cite{LT2},\ p.\ 2625) to the linear equation (\ref{01}):

\begin{lemma}
\label{lem2}Let $\sigma\left(  \theta,t\right)  $ be a smooth\ solution to the
linear\ equation $\ $%
\begin{equation}
\frac{\partial\sigma}{\partial t}\left(  \theta,t\right)  =\sigma
_{\theta\theta}\left(  \theta,t\right)  +\sigma\left(  \theta,t\right)
,\ \ \ \sigma\left(  \theta,0\right)  =\sigma_{0}\left(  \theta\right)
\label{lin}%
\end{equation}
on $S^{1}\times\lbrack0,\infty).\ $Then $\sigma\left(  \theta,t\right)  $ is
uniformly bounded on $S^{1}\times\lbrack0,\infty)$ if and only if $\int
_{0}^{2\pi}\sigma_{0}\left(  \theta\right)  d\theta=0.\ $Moreover, if
$\int_{0}^{2\pi}\sigma_{0}\left(  \theta\right)  d\theta=0,$ then
$\sigma\left(  \theta,t\right)  $ converges uniformly to the following
function%
\begin{equation}
\lim_{t\rightarrow\infty}\sigma\left(  \theta,t\right)  =\left(  \frac{1}{\pi
}\int_{0}^{2\pi}\sigma_{0}\left(  \theta\right)  \cos\theta d\theta\right)
\cos\theta+\left(  \frac{1}{\pi}\int_{0}^{2\pi}\sigma_{0}\left(
\theta\right)  \sin\theta d\theta\right)  \sin\theta,\ \ \ \theta\in S^{1}.
\label{lin1}%
\end{equation}

\end{lemma}

Applying Lemma \ref{lem2} to (\ref{01}), we obtain the uniformly convergence%
\begin{align*}
&  \lim_{t\rightarrow\infty}\left(  \frac{1}{k\left(  \theta,t\right)  }%
-\frac{L\left(  t\right)  }{2\pi}\right) \\
&  =\left[  \frac{1}{\pi}\int_{0}^{2\pi}\left(  \frac{1}{k_{0}\left(
\theta\right)  }-\frac{L\left(  0\right)  }{2\pi}\right)  \cos\theta
d\theta\right]  \cos\theta+\left[  \frac{1}{\pi}\int_{0}^{2\pi}\left(
\frac{1}{k_{0}\left(  \theta\right)  }-\frac{L\left(  0\right)  }{2\pi
}\right)  \sin\theta d\theta\right]  \sin\theta=0
\end{align*}
where the last identity is due to the identity (see \cite{GH}, p.\ 79)%
\begin{equation}
\int_{0}^{2\pi}\frac{\cos\theta}{k_{0}\left(  \theta\right)  }d\theta=\int
_{0}^{2\pi}\frac{\sin\theta}{k_{0}\left(  \theta\right)  }d\theta=0.
\end{equation}
Hence
\begin{equation}
\lim_{t\rightarrow\infty}\frac{1}{k\left(  \theta,t\right)  }=\frac{L\left(
\infty\right)  }{2\pi}=\dfrac{2A\left(  \infty\right)  }{L\left(
\infty\right)  },\ \ \ \text{uniformly in\ }\theta\in S^{1}. \label{02}%
\end{equation}
With the above $C^{0}$ convergence, together with the equation (\ref{01}),\ we
can obtain$\ C^{\infty}$ convergence of$\ 1/k\left(  \theta,t\right)  $ to the
number$\ L\left(  \infty\right)  /2\pi=2A\left(  \infty\right)  /L\left(
\infty\right)  .\ $Therefore we conclude that, if no singularity ($k=\infty
$)\ forming in finite time,\ the flow (\ref{dual1})\ converges to a round
circle with radius $L\left(  \infty\right)  /2\pi$ in $C^{\infty}$ sense.\ 

It is interesting to know that,\ with the help of Fourier series expansion
(see\ Theorem \ref{thm1}),\ the length of\ the flow (\ref{dual1})\ begins
with
\begin{equation}
L\left(  0\right)  =\int_{0}^{2\pi}u_{0}\left(  \theta\right)  d\theta
\ \ \left(  =\int_{0}^{2\pi}\left[  u_{0}\left(  \theta\right)  -a_{1}\left(
0\right)  \cos\theta-b_{1}\left(  0\right)  \sin\theta\right]  d\theta\right)
\end{equation}
and increases asymptotically to
\begin{equation}
L\left(  \infty\right)  =\sqrt{2\pi}\left(  \int_{0}^{2\pi}\left[
u_{0}\left(  \theta\right)  -a_{1}\left(  0\right)  \cos\theta-b_{1}\left(
0\right)  \sin\theta\right]  ^{2}d\theta\right)  ^{1/2}.
\end{equation}

One additional property worth mentioning\ is that the \emph{center}
(\emph{average position vector})\ of the evolving curve $\gamma_{t}\ $is
fixed.$\ $In fact, all $1/k$-type flows have this property. The\ \emph{center}%
\ of the evolving curve $\gamma_{t}$ is given by (see \cite{LT2},\ p.\ 2621)%
\begin{equation}
\frac{1}{\pi}\int_{0}^{2\pi}u\left(  \theta,t\right)  \left(  \cos\theta
,\sin\theta\right)  d\theta, \label{center}%
\end{equation}
where$\ u\left(  \theta,t\right)  $ is the support function of $\gamma_{t}%
\ $and it satisfies the evolution equation (\ref{dtu}).\ This is
\emph{independent\ of time} because it is obvious that
\begin{equation}
\frac{d}{dt}\left(  \frac{1}{\pi}\int_{0}^{2\pi}u\left(  \theta,t\right)
\left(  \cos\theta,\sin\theta\right)  d\theta\right)  =0.
\end{equation}
At this moment,\ the asymptotic behavior of the nonlocal flow (\ref{dual1}) is
well-understood (if no singularity forming in finite time).\ 

\ \ \ \ 

To end this section, we would like to say something about the Ma-Cheng's flow
(\ref{F4}).\ If there is a convex\ solution to the flow\ (\ref{F4}) on
$[0,T)$, the support function $u\left(  \theta,t\right)  $\ satisfies%
\begin{equation}
\frac{\partial u}{\partial t}\left(  \theta,t\right)  =u_{\theta\theta}\left(
\theta,t\right)  +u\left(  \theta,t\right)  -\dfrac{\int_{0}^{2\pi}\left[
u_{\theta\theta}\left(  \theta,t\right)  +u\left(  \theta,t\right)  \right]
^{2}d\theta}{\int_{0}^{2\pi}u\left(  \theta,t\right)  d\theta},\ \ \ \left(
\theta,t\right)  \in S^{1}\times\lbrack0,T) \label{MC1}%
\end{equation}
with$\ u\left(  \theta,0\right)  =u_{0}\left(  \theta\right)  ,\ u_{0}%
^{\prime\prime}\left(  \theta\right)  +u_{0}\left(  \theta\right)  >0.\ $On
the other hand, similar to equation\ (\ref{dtu}), one can use Fourier series
expansion to prove the existence of a solution $u\left(  \theta,t\right)  $ to
(\ref{MC1})\ for short time.\ This implies the existence of a convex\ solution
to the flow\ (\ref{F4}) for short time also.

We can check that if
\[
u\left(  \theta,t\right)  =\dfrac{a_{0}\left(  t\right)  }{2}+\sum
_{n=1}^{\infty}\left(  a_{n}\left(  t\right)  \cos n\theta+b_{n}\left(
t\right)  \sin n\theta\right)  ,\ \ \ u\left(  \theta,0\right)  =u_{0}\left(
\theta\right)
\]
where$\ a_{n}\left(  t\right)  =a_{n}\left(  0\right)  e^{\left(
1-n^{2}\right)  t},\ b_{n}\left(  t\right)  =b_{n}\left(  0\right)  e^{\left(
1-n^{2}\right)  t},\ $and$\ a_{0}\left(  t\right)  \ $satisfies
\[
\dfrac{a_{0}^{\prime}\left(  t\right)  }{2}=-\dfrac{\sum_{n=2}^{\infty}\left(
n^{2}-1\right)  ^{2}\left[  a_{n}^{2}\left(  t\right)  +b_{n}^{2}\left(
t\right)  \right]  }{a_{0}\left(  t\right)  }\leq0,\ \ \ \ \ a_{0}\left(
0\right)  >0
\]
then$\ u\left(  \theta,t\right)  $ will be a solution to\ (\ref{MC1}).\ We
find that%
\[
a_{0}\left(  t\right)  =\sqrt{a_{0}^{2}\left(  0\right)  -2%
{\displaystyle\sum_{n=2}^{\infty}}
\left(  n^{2}-1\right)  \left(  1-e^{2\left(  1-n^{2}\right)  t}\right)
\left(  a_{n}^{2}\left(  0\right)  +b_{n}^{2}\left(  0\right)  \right)  }.
\]
Thus (\ref{MC1})\ has a solution for short time and maintains the inequality
(\ref{uu}). Therefore the flow\ (\ref{F4}) also has a convex\ solution for
short time. Again there is a possibility that $u_{\theta\theta}\left(
\theta_{0},t_{0}\right)  +u\left(  \theta_{0},t_{0}\right)  =0\ $at finite
time $t_{0}$ for some $\theta_{0}\in S^{1},$ which is bad and we are not able
to exclude it.\ 

As $t\rightarrow\infty,$ we get%
\begin{equation}
u\left(  \theta,t\right)  \rightarrow\frac{1}{2}\sqrt{a_{0}^{2}\left(
0\right)  -2%
{\displaystyle\sum_{n=2}^{\infty}}
\left(  n^{2}-1\right)  \left(  a_{n}^{2}\left(  0\right)  +b_{n}^{2}\left(
0\right)  \right)  }+a_{1}\left(  0\right)  \cos\theta+b_{1}\left(  0\right)
\sin\theta.
\end{equation}
Note that
\[
2A\left(  0\right)  =\int_{0}^{2\pi}u_{0}\left(  \theta\right)  \left[
\left(  u_{0}\right)  _{\theta\theta}\left(  \theta\right)  +u_{0}\left(
\theta\right)  \right]  d\theta=\frac{\pi}{2}\left(  a_{0}^{2}\left(
0\right)  -2\sum_{n=2}^{\infty}\left(  n^{2}-1\right)  \left[  a_{n}%
^{2}\left(  0\right)  +b_{n}^{2}\left(  0\right)  \right]  \right)
\]
and so%
\begin{equation}
u\left(  \theta,t\right)  \rightarrow\sqrt{\frac{A\left(  0\right)  }{\pi}%
}+a_{1}\left(  0\right)  \cos\theta+b_{1}\left(  0\right)  \sin\theta
\ \ \ \text{as\ \ \ }t\rightarrow\infty.
\end{equation}
Hence the flow converges to a circle with radius $\sqrt{A\left(  0\right)
/\pi}\ $centered at $\left(  a_{1}\left(  0\right)  ,b_{1}\left(  0\right)
\right)  .\ $Its enclosed area is same as the initial area$\ A\left(
0\right)  .\ $This matches with the fact that the flow
is\ \emph{area-preserving}.$\ $Again,\ the\ center\ of the evolving curve
$\gamma_{t}$ is fixed.\ See Ma-Cheng\ \cite{MC} for more detailed discussion
of the flow.\ 

\begin{remark}
To prove asymptotic convergence\ of the flow, in p. 10, Theorem 21 of
Ma-Cheng\ \cite{MC} they quote an estimate established in Gage-Hamilton
\cite{GH}, which says that $k\left(  \theta,t\right)  r_{in}\left(  t\right)
\ $($r_{in}\left(  t\right)  $ is the inradius of $\gamma_{t}$) converges
uniformly to $1$ when the IPD $L^{2}\left(  t\right)  -4\pi A\left(  t\right)
\rightarrow0\ $(see Theorem 5.4, p. 89\ of \cite{GH})$.\ $We can not
understand why this result can be applied to their flow since these two flows
(curve shortening flow and flow (\ref{F4}))\ are very different.\ For example,
in Corollary 5.2, p. 88 of \cite{GH}\ we have \emph{Harnack-type estimate}
$k\left(  \theta,t\right)  \geq\left(  1-\varepsilon\right)  k_{\max}\left(
t\right)  $ near maximum point.\ But in \cite{MC} there is no such estimate at all.\ 
\end{remark}

\section{Determine Parallel Relation Using Mixed Isoperimetric Difference\ }

In classical differential geometry, it is known that if two simple closed
curves $\alpha,\ \beta\ $are \emph{parallel} with distance $r>0\ $%
apart\ (assume $\beta$ is an outer parallel of the fixed curve $\alpha
$),\ their length,\ enclosed area,\ and curvature are related by (see Do Carmo
\cite{D},\ p. 47)
\begin{equation}
L_{\beta}=L_{\alpha}+2\pi r,\ \ \ A_{\beta}=A_{\alpha}+rL_{\alpha}+\pi
r^{2},\ \ \ k_{\beta}\left(  s\right)  =\frac{k_{\alpha}\left(  s\right)
}{1+rk_{\alpha}\left(  s\right)  }, \label{df}%
\end{equation}
where\ $s$\ is arc\ length parameter of $\alpha.\ $As a consequence, we have
the infinitesimal identities%
\begin{equation}
\frac{dL_{\beta}}{dr}=2\pi,\ \ \ \ \ \frac{dA_{\beta}}{dr}=L_{\beta
},\ \ \ \frac{dk_{\beta}}{dr}\left(  s\right)  =-k_{\beta}^{2}\left(
s\right)  \label{P1}%
\end{equation}
and
\begin{equation}
\frac{d}{dr}\left(  L_{\beta}^{2}-4\pi A_{\beta}\right)  =0,\ \ \ \frac{d}%
{dr}\left(  \frac{L_{\beta}^{2}}{4\pi A_{\beta}}\right)  =-\frac{L_{\beta}%
}{4\pi A_{\beta}^{2}}\left(  L_{\beta}^{2}-4\pi A_{\beta}\right)  \leq0
\label{F}%
\end{equation}
for all\ $r>0\ $small (as long as the denominator of (\ref{df})\ is not zero).

By the derivative formulas in (\ref{F}),\ we clearly have:

\begin{lemma}
(\emph{parallel-invariance of the IPD})\label{lem-parallel}\ If two simple
closed curves$\ \alpha,\ \beta$ are\ \emph{parallel},\ then they have the
\emph{same} IPD, that is%
\begin{equation}
L_{\beta}^{2}-4\pi A_{\beta}=L_{\alpha}^{2}-4\pi A_{\alpha}.
\end{equation}
Moreover,\ a curve's inner\ (outer)\ parallels increase (decrease)\ its IPR.
\end{lemma}

For the convex case, there is an additional invariance, which is
\begin{equation}
\int_{0}^{2\pi}\frac{1}{k_{\beta}^{2}\left(  \theta\right)  }d\theta
-2A_{\beta}\ \ \ \left(  \text{or\ }\int_{0}^{2\pi}\frac{1}{k_{\beta}%
^{2}\left(  \theta\right)  }d\theta-\frac{L_{\beta}^{2}}{2\pi}\right)
\end{equation}
due to%
\begin{align*}
&  \int_{0}^{2\pi}\frac{1}{k_{\beta}^{2}\left(  \theta\right)  }%
d\theta-2A_{\beta}\\
&  =\int_{0}^{2\pi}\frac{1}{k_{\alpha}^{2}\left(  \theta\right)  }%
d\theta+2rL_{\alpha}+2\pi r^{2}-2\left(  A_{\alpha}+rL_{\alpha}+\pi
r^{2}\right)  =\int_{0}^{2\pi}\frac{1}{k_{\alpha}^{2}\left(  \theta\right)
}d\theta-2A_{\alpha}%
\end{align*}
where by the refined isoperimetric inequality (\ref{Tsai}) we know
\begin{equation}
\int_{0}^{2\pi}\frac{1}{k_{\alpha}^{2}\left(  \theta\right)  }d\theta
-2A_{\alpha}\geq\frac{2}{\pi}\left(  L_{\alpha}^{2}-4\pi A_{\alpha}\right)
\geq0.
\end{equation}

Combining the simple identities in (\ref{P1}) and using Gage's inequality
(\ref{Gage})\ for convex curves, we can obtain the\ interesting
\emph{monotonicity formula for convex\ parallel curves}\ (or call\emph{
}it\emph{ entropy estimate}):

\begin{lemma}
\label{lem-en}(\emph{monotonicity formula for convex\ parallel curves})\ There
holds the inequality
\begin{equation}
\frac{d}{dr}\int_{0}^{2\pi}\log\left(  k_{\beta}\left(  \theta\right)
\sqrt{\frac{A_{\beta}}{\pi}}\right)  d\theta\leq0,\ \ \ \forall\ r\in
\lbrack0,\infty). \label{mono}%
\end{equation}
Hence the integral
\begin{equation}
\int_{0}^{2\pi}\log\left(  k_{\beta}\left(  \theta\right)  \sqrt
{\frac{A_{\beta}}{\pi}}\right)  d\theta\geq0 \label{mono1}%
\end{equation}
is a \emph{decreasing} function of $r\in\lbrack0,\infty),$ which will converge
to $0$ as $r\rightarrow\infty.$
\end{lemma}

\begin{remark}
One can view Lemma \ref{lem-en} as the \emph{integration} of Gage's inequality
(under parallel evolution).\ This provides a clear explanation of Theorem 0.6
in p. 661 of Green-Osher \cite{GO},\ which has been described by them as
"physically intriguing".\ 
\end{remark}

\begin{remark}
By (\ref{mono1}) we have the \emph{entropy estimate} for a convex closed
curves\ $\gamma\ $(also see Theorem 0.2 of \cite{GO}):%
\begin{equation}
\int_{0}^{2\pi}\log\left(  k\left(  \theta\right)  \sqrt{\frac{A}{\pi}%
}\right)  d\theta\geq0,\ \ \ A=\text{enclosed area of }\gamma\label{log}%
\end{equation}
where $k\left(  \theta\right)  $ is the curvature of $\gamma$. This is already
known as a consequence of the fact that under the \emph{normalized curve
shortening flow}, the entropy is decreasing to $0\ $as $t\rightarrow\infty
\ $(see p. 10\ of the book by Zhu \cite{Z}).\ Obviously if $\gamma\ $is a
circle, then the equality holds.\ What is not so clear is the converse. If we
have equality\ in (\ref{log}), then by the above lemma we must have equalities
in both (\ref{mono})\ and (\ref{mono1}). In particular we have
\begin{equation}
\int_{0}^{2\pi}k\left(  \theta\right)  d\theta=\frac{\pi L}{A}%
,\ \ \ L=\text{length of }\gamma. \label{G}%
\end{equation}
This is the equality case of Gage's inequality, which is not discussed in
Gage's paper \cite{GA3} either.\ But according to a recent communication with
Professor Gage, he asserted that (\ref{G})\ implies $\gamma$ is a
circle.\ Hence the equality holds in (\ref{log})\ if and only if $\gamma$ is a circle.\ 
\end{remark}

\begin{remark}
Is there a proof of the entropy estimate (\ref{log}) without using a flow method?\ 
\end{remark}

%

\proof
By Gage's inequality for convex curves\ (Gage's inequality is not true for
non-convex curves) we have%
\begin{equation}
\frac{d}{dr}\int_{0}^{2\pi}\left[  \log k_{\beta}\left(  \theta\right)
+\frac{1}{2}\log\left(  \frac{A_{\beta}}{\pi}\right)  \right]  d\theta
=-\int_{0}^{2\pi}k_{\beta}\left(  \theta\right)  d\theta+\frac{\pi L_{\beta}%
}{A_{\beta}}\leq0. \label{gage}%
\end{equation}
Also note that%
\[
\lim_{r\rightarrow\infty}\int_{0}^{2\pi}\log\left(  k_{\beta}\left(
\theta\right)  \sqrt{\frac{A_{\beta}}{\pi}}\right)  d\theta=\lim
_{r\rightarrow\infty}\int_{0}^{2\pi}\log\left(  \frac{k_{\alpha}\left(
s\right)  }{1+rk_{\alpha}\left(  s\right)  }\sqrt{\frac{A_{\alpha}+rL_{\alpha
}+\pi r^{2}}{\pi}}\right)  d\theta=0.
\]
The proof is done.$%
\hfill
\square$

\ \ \ 

The converse of Lemma \ref{lem-parallel} is clearly not true.\ For two convex
curves\ $C_{1}\ $and $C_{2}$, we will show that if their IPD and \emph{mixed
IPD\ }are all the same, then they must be parallel.\ This is\ motivated by the
fact that if their IPR and \emph{mixed IPR} are all the same, then they must
be homothetic. See Lemma \ref{lem2-1} and\ Lemma \ref{lem3} below.

Recall that two simple closed curves $\alpha\left(  s\right)  ,\ \beta\left(
s\right)  :I\rightarrow\mathbb{R}^{2},$ are said to be \emph{homothetic} if
there exist some constant $\lambda>0$ and some point $\left(  a,b\right)
\in\mathbb{R}^{2}$ such that$\ \beta\left(  s\right)  =\lambda\alpha\left(
s\right)  +\left(  a,b\right)  \ $for all$\ s\in I.$ Clearly two homothetic
curves have the same IPR.$\ $This is similar to the property that two parallel
curves have the same IPD.

If two simple closed curves $\alpha,\ \beta\ $have the same IPR and IPD, then
\
\[
\left(  4\pi A_{\beta}-4\pi A_{\alpha}\right)  \left(  \frac{L_{\beta}^{2}%
}{4\pi A_{\beta}}-1\right)  =0.
\]
Thus if $\beta$ is not a circle, then$\ A_{\alpha}=A_{\beta}$ and
also$\ L_{\alpha}=L_{\beta}.\ $But if $\beta$ is a circle, so is $\alpha
.\ $Thus unless they are both circles, they must have the same length and the same\ enclosed\ area.\ 

Note that even for the convex case, two different curves $\alpha,\ \beta\ $may
have the same length and enclosed area but without other
significant\ relations at all.\ Let $u_{\alpha}\ $be the support function of
the convex\ curve $\alpha\ $with Fourier series $\ \ $
\begin{equation}
u_{\alpha}\left(  \theta\right)  =\frac{a_{0}}{2}+\sum_{n=1}^{\infty}\left(
a_{n}\cos n\theta+b_{n}\sin n\theta\right)  . \label{ua}%
\end{equation}
If we replace the coefficients $a_{n},\ b_{n}$ by $-a_{n},\ -b_{n}$ in
(\ref{ua}) for some $n$, where $n$ is large enough to maintain the convexity
condition $u_{\theta\theta}\left(  \theta\right)  +u\left(  \theta\right)
>0,$ then the new convex curve $\beta$ (with the new support function)\ will
have the same length and area as $\alpha\ $due to formulas (\ref{08-1})\ and
(\ref{08}). The curve $\beta$ is just a small perturbation of the
curve$\ \alpha.\ $

We now ask the following two interesting converse questions:

\begin{enumerate}
\item[(A).] If two curves have the same IPR,\ under what conditions are they homothetic\ ?\ 

\item[(B).] If two curves have the same IPD,\ under what conditions are they parallel\ ?
\end{enumerate}

When the curves are \emph{convex}\ (with positive curvature\ everywhere), we
can answer these questions. For general case of simple closed curves, we still
do not know the answer.

We now confine to the convex case. The advantage is that any convex curve can
be parametrized by its outward normal angle $\theta\in\left[  0,2\pi\right]  $
and its support function $u\left(  \theta\right)  \ $has domain $\theta
\in\left[  0,2\pi\right]  .\ $Moreover, $u\left(  \theta\right)  \ $behaves
well with respect to either homothetic or parallel relation. If two convex
curves $C_{1},\ C_{2}$ are homothetic, their support functions $u_{1}\left(
\theta\right)  ,\ u_{2}\left(  \theta\right)  \ $satisfy the identity$\ $%
\begin{equation}
u_{2}\left(  \theta\right)  =\lambda u_{1}\left(  \theta\right)  +a\cos
\theta+b\sin\theta,\ \ \ \forall\ \theta\in\left[  0,2\pi\right]
\end{equation}
for some constants $\lambda>0,\ a,\ b\in\mathbb{R}.\ $If they are parallel,
then $u_{2}\left(  \theta\right)  =u_{1}\left(  \theta\right)  +r$ for some
constant $r\in\mathbb{R}.\ $

It is also known that if a$\ 2\pi$-periodic function $u\left(  \theta\right)
$ satisfies the inequality$\ $%
\begin{equation}
u_{\theta\theta}\left(  \theta\right)  +u\left(  \theta\right)
>0,\ \ \ \forall\ \theta\in\left[  0,2\pi\right]  \label{uupos}%
\end{equation}
then it becomes the support function of a simple convex closed curve $C$ in
the plane (see \cite{LT2}).\ The parametrization of $C$ is given by%
\begin{equation}
X\left(  \theta\right)  =u\left(  \theta\right)  \left(  \cos\theta,\sin
\theta\right)  +u_{\theta}\left(  \theta\right)  \left(  -\sin\theta
,\cos\theta\right)  ,\ \ \ \theta\in\left[  0,2\pi\right]  . \label{U}%
\end{equation}
and its curvature is given by$\ k\left(  \theta\right)  =1/\left[
u_{\theta\theta}\left(  \theta\right)  +u\left(  \theta\right)  \right]
>0,\ \theta\in\left[  0,2\pi\right]  .$

Another nice property for the\ support functions is related to the \emph{sum}
and \emph{mixed area }of the regions enclosed by convex curves\ $C_{1}%
,\ C_{2}.\ $Let $\Omega_{1},\ \Omega_{2}$ be two strictly convex plane regions
enclosed by $C_{1},\ C_{2}.$\ The \emph{sum} (vector sum in $\mathbb{R}^{2}%
$)\ of $\Omega_{1}\ $and $\Omega_{2}\ $is$\ $defined by%
\begin{equation}
\Omega:=\Omega_{1}+\Omega_{2}=\left\{  a+b:a\in\Omega_{1}\subset\mathbb{R}%
^{2},\ b\in\Omega_{2}\subset\mathbb{R}^{2}\right\}  \subset\mathbb{R}^{2}.
\label{sum}%
\end{equation}
It is easy to check that $\Omega$ is also a convex region in $\mathbb{R}^{2}.$
Moreover, in classical convex geometry, it is proved that the boundary of
$\Omega$ is a convex\ closed\ curve with support function $u\left(
\theta\right)  \ $satisfying $u\left(  \theta\right)  =u_{1}\left(
\theta\right)  +u_{2}\left(  \theta\right)  $ for all $\theta\in\left[
0,2\pi\right]  .$\ The unique\ boundary\ point$\ $of $\Omega\ $with outward
normal angle $\theta\ $comes from the sum of the unique point on $C_{1}\ $with
outward normal angle $\theta$\ and the unique point on $C_{2}\ $with outward
normal angle $\theta,\ $and$\ $no others.

By (\ref{sum}), we can define the \emph{mixed area\ }$A\left(  \Omega
_{1},\Omega_{2}\right)  $\emph{ }of $\Omega_{1}\ $and$\ \Omega_{2},$ which is
through the identity%
\begin{equation}
A\left(  \Omega_{1}+\Omega_{2}\right)  =A\left(  \Omega_{1}\right)  +2A\left(
\Omega_{1},\Omega_{2}\right)  +A\left(  \Omega_{2}\right)  . \label{MA}%
\end{equation}
Using the support functions $u_{1}\left(  \theta\right)  ,\ u_{2}\left(
\theta\right)  ,\ $we have
\begin{equation}
A\left(  \Omega_{1}+\Omega_{2}\right)  =\frac{1}{2}\int_{0}^{2\pi}\left(
u_{1}\left(  \theta\right)  +u_{2}\left(  \theta\right)  \right)  \left[
\left(  u_{1}^{\prime\prime}\left(  \theta\right)  +u_{1}\left(
\theta\right)  \right)  +\left(  u_{2}^{\prime\prime}\left(  \theta\right)
+u_{2}\left(  \theta\right)  \right)  \right]  d\theta. \label{A12}%
\end{equation}
Hence the mixed area$\ A\left(  \Omega_{1},\Omega_{2}\right)  $ is\ given by
\begin{align}
&  2A\left(  \Omega_{1},\Omega_{2}\right) \nonumber\\
&  =\frac{1}{2}\int_{0}^{2\pi}u_{1}\left(  \theta\right)  \left(
u_{2}^{\prime\prime}\left(  \theta\right)  +u_{2}\left(  \theta\right)
\right)  d\theta+\frac{1}{2}\int_{0}^{2\pi}u_{2}\left(  \theta\right)  \left(
u_{1}^{\prime\prime}\left(  \theta\right)  +u_{1}\left(  \theta\right)
\right)  d\theta:=A_{12}+A_{21}, \label{2A12}%
\end{align}
where,\ by integration by parts, we actually have $A_{12}=A_{21}$.\ From now
on, we shall denote the mixed area $A\left(  \Omega_{1},\Omega_{2}\right)  $
by $A_{12}.\ $

Note that the values of the integrals in (\ref{A12})\ and (\ref{2A12})\ are
invariant under the$\ $transformation$\ $%
\begin{equation}
u_{1}\left(  \theta\right)  \rightarrow u_{1}\left(  \theta\right)  +a_{1}%
\cos\theta+b_{1}\sin\theta,\ \ \ \ \ u_{2}\left(  \theta\right)  \rightarrow
u_{2}\left(  \theta\right)  +a_{2}\cos\theta+b_{2}\sin\theta, \label{tr}%
\end{equation}
where$\ a_{1},\ b_{1},\ a_{2},\ b_{2}$ are constants.\ Geometrically, this
says that the areas $A\left(  \Omega_{1}+\Omega_{2}\right)  \ $and\ $A\left(
\Omega_{1},\Omega_{2}\right)  \ $are invariant under translations of
$\Omega_{1}\ $and $\Omega_{2}$ in $\mathbb{R}^{2}.\ $

With the mixed area $A_{12}$, we can define the\ \emph{mixed IPD}%
$\ $as$\ L_{1}L_{2}-4\pi A_{12},\ $where $L_{1},\ L_{2}$ are the length of
$C_{1},\ C_{2}$ respectively. Unlike the usual IPD\ $L^{2}-4\pi A\geq0$, the
mixed IPD\textsf{\ }$L_{1}L_{2}-4\pi A_{12}$ can be positive or negative\ (see
Theorem \ref{thm2} below).\ In particular, if one of $C_{1},\ C_{2}\ $is a
circle,$\ $say$\ C_{1}$ is a circle with radius $r$, then $A_{12}%
=A_{21}=rL_{2}/2$ and\ the mixed IPD$\ $disappears with$\ L_{1}L_{2}-4\pi
A_{12}=2\pi rL_{2}-2\pi rL_{2}=0.\ $The mixed IPD is related to the IPD of
$\Omega_{1}+\Omega_{2}$ via the identity%
\begin{equation}
L^{2}-4\pi A\ \text{(for }\Omega_{1}+\Omega_{2}\text{)}=\left(  L_{1}^{2}-4\pi
A_{1}\right)  +\left(  L_{2}^{2}-4\pi A_{2}\right)  +2\left(  L_{1}L_{2}-4\pi
A_{12}\right)  . \label{LLL}%
\end{equation}

By (\ref{LLL}) and Lemma \ref{lem-parallel}, the mixed IPD $L_{1}L_{2}-4\pi
A_{12}\ $is\emph{ }clearly\emph{ invariant} under parallel relations also.
That is:

\begin{lemma}
\label{lem-1}(\emph{parallel-invariance of} \emph{the} \emph{mixed IPD})\ If
we replace $C_{1}$ by a parallel curve $\tilde{C}_{1}$ and $C_{2}$ by another
parallel curve $\tilde{C}_{2},\ $then the mixed IPD$\ $for the pair$\ \tilde
{C}_{1},$ $\tilde{C}_{2}$ is the same as that for the pair$\ C_{1},$ $C_{2}$.
\end{lemma}

The most important property for mixed area is the following well-known
\emph{Minkowski\ mixed area\ inequality\ }(see the encyclopedic\ book by
R.\ Schneider \cite{S}):
\begin{equation}
\sqrt{A_{1}A_{2}}\leq A_{12}, \label{AA}%
\end{equation}
where the\ equality holds\ if and only if $C_{1}$ and $C_{2}$ are
\emph{homothetic}.

With the above theorem, our answer to question (A) comes immediately:

\begin{lemma}
(\emph{characterization of} \emph{homothetic relation})\label{lem2-1}\ Two
convex closed curves $C_{1},\ C_{2}\ $are \emph{homothetic }if and only if
\begin{equation}
\frac{L_{1}^{2}}{4\pi A_{1}}=\frac{L_{2}^{2}}{4\pi A_{2}}=\frac{L_{1}L_{2}%
}{4\pi A_{12}}. \label{3-same}%
\end{equation}
That is, all IPR, including the mixed one, are the same.\ 
\end{lemma}

\begin{remark}
Note that the \emph{mixed IPR}$\ L_{1}L_{2}/4\pi A_{12}$ is invariant under
dilations or translations of $C_{1}\ $and$\ C_{2}$.\ Compare with the
parallel-invariance property in Lemma \ref{lem-1}.
\end{remark}

%

\proof
If$\ C_{1},\ C_{2}\ $are homothetic, we clearly have$\ L_{1}^{2}/4\pi
A_{1}=L_{2}^{2}/4\pi A_{2}\ $(call this value $\lambda$),$\ $which
gives$\ L_{1}=\sqrt{4\pi\lambda A_{1}},\ L_{2}=\sqrt{4\pi\lambda A_{2}}.$ By
Minkowski inequality, we also have $A_{12}=\sqrt{A_{1}A_{2}}.$ Hence%
\[
\frac{L_{1}L_{2}}{4\pi A_{12}}=\frac{\sqrt{4\pi\lambda A_{1}}\sqrt{4\pi\lambda
A_{2}}}{4\pi\sqrt{A_{1}A_{2}}}=\lambda.
\]
Conversely if (\ref{3-same})\ holds\ with common value$\ \lambda$,
then$\ L_{1}=\sqrt{4\pi A_{1}\lambda},\ L_{2}=\sqrt{4\pi A_{2}\lambda},$ and%
\[
4\pi A_{12}=\frac{L_{1}L_{2}}{\lambda}=\frac{\sqrt{4\pi A_{1}\lambda}%
\sqrt{4\pi A_{2}\lambda}}{\lambda}=4\pi\sqrt{A_{1}A_{2}},
\]
which gives the equality $A_{12}=\sqrt{A_{1}A_{2}}.$ By Minkowski theorem they
are homothetic. The proof is done.$%
\hfill
\square$\ 

\ \ \ 

To answer question (B), one can not rely on the Minkowski inequality because
it does not have the right form.\ Instead we can use the following, which
concerns the IPD:

\begin{theorem}
\label{thm2}(\emph{mixed} \emph{IPD inequality})\ For any two$\ $convex closed
curves $C_{1},\ C_{2}$ in $\mathbb{R}^{2}\ $with\ support functions\ $u_{1}%
\left(  \theta\right)  ,$ $u_{2}\left(  \theta\right)  $, there holds the
estimate%
\begin{equation}
-\sqrt{I_{1}}\sqrt{I_{2}}\leq L_{1}L_{2}-4\pi A_{12}\leq\sqrt{I_{1}}%
\sqrt{I_{2}}\label{LL}%
\end{equation}
where $I_{1}=L_{1}^{2}-4\pi A_{1},\ I_{2}=L_{2}^{2}-4\pi A_{2},\ $%
and$\ A_{12}$ is the mixed area of $C_{1}\ $and $C_{2}.\ $The equality holds
in the lower bound estimate if and only if
\begin{equation}
\frac{\sqrt{I_{2}}}{k_{1}\left(  \theta\right)  }+\frac{\sqrt{I_{1}}}%
{k_{2}\left(  \theta\right)  }=c,\ \ \ \forall\ \theta\in\left[
0,2\pi\right]  \label{LL1}%
\end{equation}
for some constant\ $c>0.\ $Also, the equality holds in the upper bound
estimate if and only if$\ $
\begin{equation}
\frac{\sqrt{I_{2}}}{k_{1}\left(  \theta\right)  }-\frac{\sqrt{I_{1}}}%
{k_{2}\left(  \theta\right)  }=c,\ \ \ \forall\ \theta\in\left[
0,2\pi\right]  \label{LL2}%
\end{equation}
for some constant\ $c.\ $Here $k_{1}\left(  \theta\right)  ,\ k_{2}\left(
\theta\right)  $ are the curvature of $C_{1},\ C_{2}\ $respectively.
\end{theorem}

\begin{remark}
\label{Sch}We actually proved inequality (\ref{LL}) by ourselves (just a
simple observation), but later\ on Professor\ Schneider kindly told us that it
had been proved\ by Favard in 1930\ and also reappeared in (4),\ p. 105 of the
book\ Bonnesen-Fenchel\ \cite{BF}. A higher-dimensional version of
(\ref{LL})\ is inequality (6.4.9),\ p. 335 of his book \cite{S}.\ However,
since we do not see the equality results (\ref{LL1})\ and (\ref{LL2})\ in p.
105 of \cite{BF} or \cite{S}, we still give a proof of Theorem \ref{thm2}
because it takes only a few lines.
\end{remark}

%

\proof
Without loss of generality, we may assume $I_{1}>0,\ I_{2}>0$, otherwise we
would have $L_{1}L_{2}-4\pi A_{12}=0$ and (\ref{LL}), (\ref{LL1}), (\ref{LL2})
all hold.\ Let $u_{1}\left(  \theta\right)  ,\ u_{2}\left(  \theta\right)  $
be the support functions of $C_{1},\ C_{2}\ $respectively.\ If $C_{\alpha
\beta}$ is a convex closed curve with support function $u_{\alpha\beta}\left(
\theta\right)  $ given by the linear combination%
\[
u_{\alpha\beta}\left(  \theta\right)  =\alpha u_{1}\left(  \theta\right)
+\beta u_{2}\left(  \theta\right)  ,\ \ \ \alpha,\ \beta\ \text{are
non-zero\ constant}%
\]
then the IPD of $C_{\alpha\beta}\ $is given by%
\begin{equation}
L_{\alpha\beta}^{2}-4\pi A_{\alpha\beta}=\alpha^{2}I_{1}+\beta^{2}%
I_{2}+2\alpha\beta\left(  L_{1}L_{2}-4\pi A_{12}\right)  \geq0.
\end{equation}
Thus if we choose $\alpha=\sqrt{I_{2}},\ \beta=\sqrt{I_{1}}\ $(this is the
optimal choice),$\ $we would have%
\[
L_{\alpha\beta}^{2}-4\pi A_{\alpha\beta}=2I_{1}I_{2}+2\sqrt{I_{1}I_{2}}\left(
L_{1}L_{2}-4\pi A_{12}\right)  \geq0,
\]
which gives\ the lower bound and the equality holds if and only if
$C_{\alpha\beta}$ is a circle with constant curvature. Hence (\ref{LL}) follows.

On the other hand, if we choose$\ \alpha=\sqrt{I_{2}},\ \beta=-\sqrt{I_{1}},$
then$\ u_{\alpha\beta}\left(  \theta\right)  \ $may not satisfy $u_{\alpha
\beta}^{\prime\prime}\left(  \theta\right)  +u_{\alpha\beta}\left(
\theta\right)  >0.\ $But we can modify it by considering%
\[
u_{\alpha\beta}\left(  \theta\right)  =\sqrt{I_{2}}u_{1}\left(  \theta\right)
-\sqrt{I_{1}}u_{2}\left(  \theta\right)  +c
\]
for some large constant $c>0.\ $Now $u_{\alpha\beta}\left(  \theta\right)  $
is the support function of some convex closed curve $C_{\alpha\beta}\ $with
(adding $c$ has no effect in $L_{\alpha\beta}^{2}-4\pi A_{\alpha\beta}$)
\[
L_{\alpha\beta}^{2}-4\pi A_{\alpha\beta}=2I_{1}I_{2}-2\sqrt{I_{1}I_{2}}\left(
L_{1}L_{2}-4\pi A_{12}\right)  \geq0,
\]
which gives the upper bound.\ The equality holds if and only if $C_{\alpha
\beta}\ $is a circle and we have (\ref{LL2}).$%
\hfill
\square$

\ \ \ \ \ 

Motivated by Lemma \ref{lem2-1} and with the help of Theorem \ref{thm2}, our
answer to question (B) is given by:\ 

\begin{lemma}
(\emph{characterization of} \emph{parallel relation})\label{lem3}\ Two convex
curves $C_{1},\ C_{2}\ $are \emph{parallel} (up to a translation)\emph{ }if
and only if
\begin{equation}
L_{1}^{2}-4\pi A_{1}=L_{2}^{2}-4\pi A_{2}=L_{1}L_{2}-4\pi A_{12}. \label{3-d}%
\end{equation}
That is, all IPD, including the mixed one, are the same.\ 
\end{lemma}

%

\proof
If $C_{1},\ C_{2}\ $are parallel, then $u_{2}\left(  \theta\right)
=u_{1}\left(  \theta\right)  +r+a\cos\theta+b\sin\theta\ $for some constants
$r,\ a,\ b\ $and then $L_{2}=L_{1}+2\pi r.\ $Hence%
\[
A_{12}=\frac{1}{2}\int_{0}^{2\pi}u_{2}\left(  \theta\right)  \left(
u_{1}^{\prime\prime}\left(  \theta\right)  +u_{1}\left(  \theta\right)
\right)  d\theta=A_{1}+\frac{1}{2}rL_{1}%
\]
and$\ $%
\[
L_{1}L_{2}-4\pi A_{12}=L_{1}\left(  L_{1}+2\pi r\right)  -4\pi\left(
A_{1}+\frac{1}{2}rL_{1}\right)  =L_{1}^{2}-4\pi A_{1}.
\]
Hence\ (\ref{3-d})\ holds.

Conversely, if (\ref{3-d})\ holds, then by (\ref{LL2}) we have%
\[
U^{\prime\prime}\left(  \theta\right)  +U\left(  \theta\right)
=c,\ \ U\left(  \theta\right)  :=u_{2}\left(  \theta\right)  -u_{1}\left(
\theta\right)  ,\ \ \ \forall\ \theta\in\left[  0,2\pi\right]
\]
for some constant $c.\ $If $D>0$ is a large constant, the function $V\left(
\theta\right)  =U\left(  \theta\right)  +D\ $will be the support function of
some circle of $\mathbb{R}^{2}.$ Hence%
\[
V\left(  \theta\right)  =\rho+a\cos\theta+b\sin\theta
\]
for some constants $a,\ b,\ \rho\in\mathbb{R},$ $\rho>0.\ $As a result, we get%
\[
u_{2}\left(  \theta\right)  -u_{1}\left(  \theta\right)  =r+a\cos\theta
+b\sin\theta
\]
for some constants $a,\ b,\ r\in\mathbb{R}$ and$\ C_{1},\ C_{2}\ $are
parallel.\ The proof is done.$%
\hfill
\square\ \ $

\begin{remark}
Although\ Lemma \ref{lem-parallel} is valid for all simple closed curves, it
is not clear how to generalize Lemma \ref{lem2-1} and Lemma \ref{lem3} to the
non-convex case. For non-convex sets $\Omega_{1},\ \Omega_{2}\ $in
$\mathbb{R}^{2}\ $(bounded by simple closed curves $C_{1},\ C_{2}$)$,\ $one
can still use identity (\ref{MA}) to define their mixed area$\ A\left(
\Omega_{1},\Omega_{2}\right)  =A_{12}$. Hence we can still talk about their
mixed IPR and IPD.\ However, it is not clear whether we have good results
similar to the Minkowski inequality and Theorem \ref{thm2}.\ 
\end{remark}

In higher-dimensional space, say $\mathbb{R}^{3},\ $parallel surfaces
$S_{1},\ S_{2}\ $\emph{do not have}, in general, the same IPD any more. The
classical isoperimetric inequality for a compact connected closed surface
$S\ $is%
\begin{equation}
A^{3}\left(  S\right)  \geq36\pi V^{2}\left(  \Omega\right)  \label{ISOP}%
\end{equation}
where $A\left(  S\right)  \ $is the surface area of $S\ $and$\ \Omega$ is the
domain enclosed by it with volume$\ V\left(  \Omega\right)  .\ $Moreover, the
equality holds if and only if $S$ is a sphere.\ In view of (\ref{ISOP}), the
IPD quantity for a space surface is $A^{3}-36\pi V^{2}.$

If we assume that parallel surfaces have the same\ IPD, we would have the
infinitesimal identity
\begin{equation}
\frac{d}{dr}\left[  A^{3}\left(  S_{r}\right)  -36\pi V^{2}\left(  \Omega
_{r}\right)  \right]  =0,\ \ \ S_{0}=S,
\end{equation}
where $S_{r}$ is parallel to $S$ (for small $r$) with enclosed domain
$\Omega_{r}.\ $After computation, we would get (evaluated at $r=0$)
\[
6A^{2}\left(  S\right)  \int_{S}H\left(  p\right)  dp-72\pi V\left(
\Omega\right)  A\left(  S\right)  =0,
\]
where $H$ is the \emph{mean curvature} of $S.\ $The above is same as%
\begin{equation}
A\left(  S\right)  \int_{S}H\left(  p\right)  dp=12\pi V\left(  \Omega\right)
. \label{022}%
\end{equation}
However, we know that (\ref{022})\ does not hold for a general compact closed
surface $S.$

\bigskip\ 

\ \ \ \ \ 

\textbf{Acknowledgments.\ \ }\ \ This research is supported by\ NSC (National
Science Council)\ of Taiwan, under grant number 96-2115-M-007-010-MY3. We also
appreciate the support of NCTS\ (National Center of Theoretical Sciences)\ of
Taiwan.\ The second author is grateful to the hospitality of Professor
Shengliang Pan when he was visiting East China Normal University.\ Finally we
thank Professors M.\ Gage and R.\ Schneider for some e-mail communications.\ 

\ \ \ 

\ \ \ \ \ \

\ \ \ \ 

\ \ \ \ \ 

Yu-Chu Lin

Department of Mathematics

National Tsing Hua University

Hsinchu 30013,\ TAIWAN

E-mail:\ \textit{yclin@math.nthu.edu.tw}

\ \ \ \ 

\ \ \ \ 

Dong-Ho Tsai\ 

Department of Mathematics

National Tsing Hua University

Hsinchu 30013,\ TAIWAN

E-mail:\ \textit{dhtsai@math.nthu.edu.tw}

\bigskip\

\end{document}